\documentclass[10pt]{article}

\RequirePackage{amsmath}
\RequirePackage{amsfonts}
\RequirePackage{amssymb}
\RequirePackage{hyperref}
\RequirePackage{multirow}
\RequirePackage{algorithm}
\RequirePackage{algpseudocode}
\RequirePackage{graphicx}
\RequirePackage{xcolor}
\RequirePackage{amsthm}
\RequirePackage{amsmath}
\RequirePackage{amsfonts}
\RequirePackage{amssymb}
\RequirePackage{hyperref}
\RequirePackage{multirow}
\RequirePackage{algorithm}
\RequirePackage{algpseudocode}
\RequirePackage{graphicx}
\RequirePackage{xcolor}
\allowdisplaybreaks

\usepackage{dpr_fec,ifthen,dsfont}
\hypersetup{colorlinks}

\usepackage{isetechreport}
\def\reportyear{22T}
\def\reportno{006}
\def\revisionno{0}
\def\originaldate{\today}
\def\revisiondate{\today}
\coraltrue
\cvcrfalse

\begin{document}

\title{Worst-Case Complexity of TRACE with Inexact Subproblem Solutions for Nonconvex Smooth Optimization}

\author{Frank E.~Curtis\thanks{E-mail: frank.e.curtis@lehigh.edu}}
\author{Qi Wang\thanks{E-mail: qiw420@lehigh.edu}}
\affil{Department of Industrial and Systems Engineering, Lehigh University}
\titlepage

\maketitle

\begin{abstract}
  An algorithm for solving nonconvex smooth optimization problems is proposed, analyzed, and tested.  The algorithm is an extension of the Trust Region Algorithm with Contractions and Expansions (TRACE) [Math.~Prog.~162(1):132, 2017].  In particular, the extension allows the algorithm to use inexact solutions of the arising subproblems, which is an important feature for solving large-scale problems.  Inexactness is allowed in a manner such that the optimal iteration complexity of $\Ocal(\epsilon^{-3/2})$ for attaining an $\epsilon$-approximate first-order stationary point is maintained while the worst-case complexity in terms of Hessian-vector products may be significantly improved as compared to the original TRACE.  Numerical experiments show the benefits of allowing inexact subproblem solutions and that the algorithm compares favorably to a state-of-the-art technique.
\end{abstract}

\newcommand{\ARC}{\textsc{arc}}
\newcommand{\FDS}{\textsc{fds}}
\newcommand{\gLip}{g_{{\rm Lip}}}
\newcommand{\HLip}{H_{{\rm Lip}}}
\newcommand{\HLoc}{H_{{\rm Loc}}}
\newcommand{\ITRACE}{\textsc{i-trace}}
\newcommand{\TLTR}{\textsc{tltr}}
\newcommand{\TRACE}{\textsc{trace}}

\section{Introduction}\label{sec.introduction}

There are a variety of algorithmic methodologies for solving nonconvex smooth optimization problems that offer state-of-the-art performance when solving broad subclasses of such problems.  Among these, a few offer a worst-case performance guarantee for achieving approximate first-order stationarity that is optimal with respect to a class of second-order-derivative-based methods for minimizing sufficiently smooth objective functions.  These include certain cubic-regularization, quadratic-regularization, line-search, and trust-region methods; see Section~\ref{sec.literature_review}.

In this paper, we propose, analyze, and provide the results of numerical experiments with an extended version of the \emph{Trust Region Algorithm with Contractions and Expansions} (\TRACE) from \cite{CurtRobiSama17}, which was the first trust-region method to attain the aforementioned optimal iteration complexity guarantees.  In particular, the algorithm that we propose overcomes the main deficiency of \TRACE, namely, that \TRACE{} requires exact solutions of the arising trust-region subproblems, which is impractical in large-scale settings.  Our algorithm overcomes this deficiency by employing an iterative linear algebra technique---specifically, a Krylov subspace method---for solving the arising subproblems and allowing the ``outer'' algorithm for solving the original problem to use inexact solutions from the ``inner'' algorithm for solving the trust-region subproblems.  This represents a response to the conjecture from \cite{CurtRobiSama17}, which states: ``We expect that such variations of our algorithm can be designed that maintain our global convergence guarantees... our worst-case complexity bounds and local convergence guarantees.''  In fact, our proposed algorithm goes beyond this conjecture.  Not only do we show that our approach maintains the global convergence, local convergence, and worst-case complexity guarantees of \TRACE; we also show that our proposed enhancement of \TRACE{} achieves strong worst-case complexity properties in terms of the overall required number of Hessian-vector products, which are the most expensive operations required when solving many large-scale problems.

Our numerical experiments show that our algorithm offers practical benefits over another optimal-worst-case-complexity method for solving nonconvex smooth optimization problems.  We also demonstrate that our approach offers the computational flexibility in terms of the trade-offs between derivative evaluations and Hessian-vector products that should be expected of any such method that allows inexact subproblem solutions.  In particular, with more exact subproblem solutions, our algorithm often requires fewer derivative evaluations at the expense of more Hessian-vector products, whereas with more inexact subproblem solutions, it often requires fewer Hessian-vector products at the expense of more derivative evaluations.  This allows any user of our algorithm to tailor its use depending on the relative costs of these operations.

\subsection{Notation, Problem Formulation, and Assumptions}\label{sec.notation}

We use $\R{}$ to denote the set of real numbers, $\R{}_{\geq0}$ (resp.,~$\R{}_{>0}$) to denote the set of nonnegative (resp.,~positive) real numbers, $\R{n}$ to denote the set of $n$-dimensional real vectors, $\R{m \times n}$ to denote the set of $m$-by-$n$-dimensional real matrices, $\mathbb{S}^n \subset \R{n \times n}$ to denote the set of $n$-by-$n$-dimensional real symmetric matrices, and $\N{}$ to denote the set of nonnegative integers.  We use $I$ to denote the identity matrix and use $e_j$ for $j \in \N{} \setminus \{0\}$ to denote the $j$th column of the identity matrix, where in each case the dimension of the object is determined by the context in which it appears.  We use the function $|\cdot|$ to take the absolute value of a real number and the function $\|\cdot\|$ to take the 2-norm of a real vector or to take the induced 2-norm of a real matrix.  Given real numbers $a$ and $b$, we use $a \perp b$ to mean that $ab = 0$.  Given $H \in \mathbb{S}^n$, we use $H \succ 0$ (resp.,~$H \succeq 0$) to indicate that $H$ is positive definite (resp.,~semidefinite).

Given functions $\phi : \R{} \to \R{}_{\geq0}$ and $\varphi : \R{} \to \R{}_{\geq0}$, the expression $\phi(\cdot) = \Ocal(\varphi(\cdot))$ means that there exists $c \in \R{}_{>0}$ such that $\phi(\cdot) \leq c\varphi(\cdot)$.  Similarly, given positive real number sequences $\{\phi_k\}$ and $\{\varphi_k\}$, the expression $\phi_k = \Ocal(\varphi_k)$ means that there exists $c \in \R{}_{>0}$ such that $\phi_k \leq c \varphi_k$ for all $k \in \N{}$.  If, in addition, the sequences have the property that $\{\phi_k/\varphi_k\} \to 0$, then one writes $\phi_k = o(\varphi_k)$.

Our problem of interest is the minimization problem
\begin{equation}\label{prob.f}
  \min_{x\in \R{n}} f(x),
\end{equation}
where $f : \R{n} \to \R{}$ satisfies Assumption~\ref{ass.f}, stated below.  The gradient and Hessian functions for $f$ are denoted by $g := \nabla f : \R{n} \to \R{}$ and $H := \nabla^2 f : \R{n} \to \mathbb{S}^n$, respectively.  Given the $k$th iterate in an algorithm for solving \eqref{prob.f}, call it $x_k \in \R{n}$, we define $f_k := f(x_k)$, $g_k := g(x_k) \equiv \nabla f(x_k)$, and $H_k := H(x_k) \equiv \nabla^2f(x_k)$.  We also apply a subscript to refer to other quantities corresponding to the $k$th iteration; e.g., the iterate displacement (i.e., step) is denoted as $s_k \in \R{n}$.  Our algorithm involves subroutines that have their own iteration indices, and we use additional subscripts to keep track of quantities associated with the inner iterations of these subroutines.

Assumption~\ref{ass.f}, below, is made throughout the paper.  For \TRACE{} in \cite{CurtRobiSama17}, a guarantee of convergence from remote starting points is first proved under a weaker assumption---namely, without the assumption of Lipschitz continuity of the Hessian function---prior to worst-case iteration complexity guarantees being proved under an assumption on par with Assumption~\ref{ass.f}.  We claim that the same could be done for the algorithm proposed in this paper, but for the sake of brevity we jump immediately to Assumption~\ref{ass.f} in order to prove worst-case complexity properties.  As stated later, these properties, in turn, ensure convergence from remote starting points.

Assumption~\ref{ass.f} refers, for each index tuple $(k,j,l) \in \N{} \times \N{} \times \N{}$ generated by our algorithm, to the iterate~$x_k \in \R{n}$ and trial step $Q_{k,j}t_{k,j,l} \in \R{n}$.  The fact that the assumption refers to these algorithmic quantities should not be seen as a deficiency of our analysis.  After all, our algorithm guarantees monotonic nonincrease of the objective values, meaning that $\{x_k\}$ is contained in the sublevel set for $f$ with respect to the initial value $f(x_0)$, i.e., $\Lcal_{f,0} := \{x \in \R{n} : f(x) \leq f(x_0)\}$.  In addition, with each accepted step, the algorithm requires decrease in the objective that is proportional to the cubed norm of the step, so with an objective that is bounded below, it is reasonable to assume that the accepted steps are bounded in norm, which in turn means that it is reasonable to assume (by construction of our algorithm) that all trial steps are bounded in norm.  Assuming that is the case, the open convex set mentioned in the assumption would itself be contained in the Minkowski sum of $\Lcal_{f,0}$ and a bounded set, in light of which the assumption is standard for smooth optimization.

\bassumption\label{ass.f}
  The function $f : \R{n} \to \R{}$ is twice-continuously differentiable and bounded below by a real number $f_{\inf} \in \R{}$ over~$\R{n}$.  In addition, the gradient function $g : \R{n} \to \R{n}$ and Hessian function $H : \R{n} \to \mathbb{S}^n$ are each Lipschitz continuous with Lipschitz constants denoted by $\gLip \in \R{}_{>0}$ and $\HLip \in \R{}_{>0}$, respectively, in an open convex set containing $x_k$ and $x_k + Q_{k,j}t_{k,j,l}$ for all generated $(k,j,l) \in \N{} \times \N{} \times \N{}$.
\eassumption

\noindent
A consequence of the Lipschitz continuity of the gradient in Assumption~\ref{ass.f} is that the Hessian matrix $H_k$ is bounded in norm for all $k \in \N{}$ in the sense that there exists $H_{\max} \in \R{}_{>0}$ such that $\|H_k\|\leq H_{\max}$ for all $k \in \N{}$.

\subsection{Literature review}\label{sec.literature_review}

Our focus in this paper is on worst-case complexity bounds for a second-order-derivative-based algorithm to reach an iterate $x_k$ that is $\epsilon$-approximate first-order stationary (a property to which we refer as $\epsilon$-stationary throughout the paper) with respect to \eqref{prob.f} in the sense that
\bequation\label{eq.eps_stationary}
  \|\nabla f(x_k)\| \leq \epsilon.
\eequation
Some research articles have also considered worst-case complexity bounds for achieving second- or even higher-order stationarity, but since our focus in this paper is on large-scale settings in which such guarantees are impractical to require, we consider such complexity bounds outside of our scope.

Trust-region methods that employ a traditional updating scheme for the trust-region radius based on actual-to-quadratic-model-predicted-reduction ratios (see, e.g., \cite{ConGT00,NoceWrig06}) are known to have a worst-case \emph{iteration} complexity (and, correspondingly, \emph{function-} and \emph{derivative-evaluation} complexity) of $\Ocal(\epsilon^{-2})$ for achieving $\epsilon$-stationarity \cite{CartGoulToin10,CurtLubbRobi18,GrapYuanYuan15}.  Importantly, this complexity is known to be tight for both first- and second-order variants of such methods~\cite{CartGoulToin10}.  It was first shown by Nesterov and Polyak that cubic regularization of a second-order method---an idea that has appeared as far back as Griewank~\cite{Grie81}---can achieve an improved iteration complexity of $\Ocal(\epsilon^{-3/2})$; see \cite{NestPoly06}.  This complexity for achieving $\epsilon$-stationarity is now known to be optimal with respect to a class of second-order-derivative-based methods for minimizing sufficiently smooth objectives~\cite{CartGoulToin10}.  After \cite{NestPoly06}, practical variants of the idea subsequently appeared in papers by Cartis, Gould, and Toint \cite{CartGoulToin11,CartGoulToin11b} in the form of the adaptive-regularisation-using-cubics (\ARC) method.  Since that time, a few other schemes have been developed that build upon cubic- or even quadratic-regularization techniques; see, e.g., \cite{BirgMart17,CurtRobiSama18,Duss17,DussOrba15}.  These ideas have also been extended to high-order regularization of higher-order methods in order to achieve improved complexity bounds; see, e.g., the work by Birgin et al.~in~\cite{BirgGardMartSantToin17}, where it is shown that a $p$th-order method with $(p+1)$th-order regularization can achieve an iteration complexity of $\Ocal(\epsilon^{-(p+1)/p})$.

Our work in this paper is motivated by empirical observations (see, e.g., Section~\ref{sec.numerical}) that while cubic-regularization and related methods can offer strong complexity properties, they can disappoint in practice compared to traditional trust-region algorithms.  (One can design cubic-regularization methods that are competitive with trust-region methods in practice, although this requires more complicated updating schemes for the regularization parameter than are analyzed in the papers cited above~\cite{GoulPorcToin12}.)  Therefore, we contend that it is of interest to explore trust-region methods that do not deviate too much from traditional schemes, yet offer optimal worst-case complexity properties.  The first such trust-region method to achieve $\Ocal(\epsilon^{-3/2})$ iteration complexity with respect to $\epsilon$-stationarity was \TRACE{} \cite{CurtRobiSama17}.  Another approach, which tries to adhere closely to the popular combination of a trust-region method that employs the linear conjugate gradient (CG) method to solve the arising subproblems, is that in \cite{curtis2021trust}; see also the prior line-search method proposed by Royer and Wright in \cite{royer2018complexity}.

As previously mentioned, our work in this paper is motivated by the goal to improve the \emph{computational} complexity of \TRACE{} in large-scale settings where matrix factorizations and/or Hessian-vector products can dominate the computational expense.  The use of iterative linear algebra techniques to exploit potentially inexact subproblem solutions has been a topic of research for decades.  Traditional line-search and trust-region methods that use CG to solve the arising subproblems (approximately) have been studied and implemented widely \cite{Ste83,Toi81}.  It is well known that with sufficiently exact subproblem solutions, such algorithms can attain the superlinear or quadratic rates of local convergence of Newton's method \cite{DemES82}.  Similar guarantees have also been shown for \ARC{} \cite{CartGoulToin11}; see also \cite{BirgGardMartSantToin17} for the use of inexact subproblem solutions in higher-order regularization schemes.  Particularly in the case of trust-region methods when one aims to be able to solve the subproblems to arbitrary accuracy, the use of the Lanczos method has been well studied \cite{GouLRT99}, for which it is known that if the solution of a trust-region subproblem in $n$ variables lies on the boundary of the trust-region radius, then the subproblem is equivalent to an extremal eigenvalue problem of a matrix of size $2n$ \cite{adachi2017solving}.  Convergence of the Lanczos method for estimating eigenvalues has been analyzed in \cite{kuczynski1992estimating,kuczynski1994probabilistic}, the results of which have been used in the analysis of various optimization algorithms; see, e.g., \cite{carmon2018accelerated,curtis2021trust,royer2020newton,royer2018complexity}.  Complexity guarantees for the Lanczos method specifically for trust-region methods has been studied in \cite{carmon2018analysis,gould2020error,JiaW21,ZhaSL17}.  The results in \cite{gould2020error} play an important role in this paper.

\subsection{Contributions}\label{sec.contributions}

The work in this paper builds on the ideas and analyses provided in the aforementioned literature, but offers a unique contribution since we provide the first inexact variant of \TRACE, a method that we call \ITRACE, that offers iteration and gradient evaluation complexity bounds that match those of \TRACE.  We also use results about the complexity of the Lanczos algorithm to show that \ITRACE{} offers state-of-the-art complexity in terms of Hessian-vector products when solving large-scale problem instances.  Our theoretical analyses are backed by empirical evidence showing that our proposed \ITRACE{} method offers computational flexibility beyond that offered by \TRACE{} and compares favorably against an implementation of \ARC{} that also allows inexact subproblem solutions.  We attribute this behavior to the fact that \ITRACE{} adheres closely to a traditional trust-region strategy, where to achieve optimal iteration complexity it adaptively uses a combination of explicit and implicit regularization of the Hessian matrices in the arising subproblems.

\subsection{Organization}\label{sec.organization}

Section~\ref{sec.algorithm} contains a description of our algorithm and its associated subroutines.  (Our description involves well-known characterizations and properties of Krylov subspace, specifically Lanczos-based, iterative methods for solving subproblems arising in optimization algorithms, for which we refer the reader to \cite{ConGT00a,GouLRT99,ZhaSL17} and other provided references.)  Section~\ref{sec.analysis} contains our convergence and worst-case complexity analyses of the algorithm.  The results of numerical experiments are provided in Section~\ref{sec.numerical} and concluding remarks are provided in Section~\ref{sec.conclusion}.

\section{Algorithm Description}\label{sec.algorithm}

Each iteration of \ITRACE{} involves the minimization of a second-order Taylor series model of $f$ at the current iterate within a trust region; specifically, in iteration $k \in \N{}$, the model $m_k : \R{n} \to \R{}$ is defined by
\bequationNN
  m_k(s) = f_k + g_k^Ts + \thalf s^TH_ks.
\eequationNN
Building on \TRACE, the trust region is defined either explicitly through a trust-region radius $\delta \in \R{}_{>0}$ and a trust-region constraint of the form $\|s\| \leq \delta$ or implicitly through a regularization parameter $\lambda \in \R{}_{>0}$ and a regularization term $\thalf \lambda \|s\|^2$, where $\lambda$ is sufficiently large such that $H_k + \lambda I \succ 0$, which in turn means that the regularized model $m_k(\cdot) + \thalf \lambda \|\cdot\|^2 = f_k + g_k^T(\cdot) + \thalf (\cdot)^T(H_k + \lambda I)(\cdot)$ is strongly convex.

However, unlike \TRACE, the main idea behind \ITRACE{} is to allow an approximate subproblem solution to be considered acceptable.  Specifically, in each iteration $k \in \N{}$, the algorithm might only consider, for some $j \in \{0,\dots,n-1\}$, the solution of a subproblem over the $j$th-order Krylov subspace defined by $g_k$ and $H_k$, namely,
\bequationNN
  \Kcal_{k,j} := \linspan\{g_k,H_kg_k,\dots,H_k^jg_k\}.
\eequationNN
Using the Lanczos process, \ITRACE{} iteratively constructs orthonormal bases for such subspaces for increasing $j$, as needed.  Let such a basis be given by
\bequationNN
  Q_{k,j} := \bbmatrix q_{k,0} & q_{k,1} & \cdots & q_{k,j} \ebmatrix \in \R{n \times (j+1)}.
\eequationNN
With this basis constructed using Lanczos, one finds for any $(k,j) \in \N{} \times \{0,\dots,n-1\}$ that there exists tridiagonal $T_{k,j} \in \R{(j+1) \times (j+1)}$ such that with $\gamma_{k,0} \gets \|g_k\|$ one has
\bequationNN
  T_{k,j} = Q_{k,j}^TH_kQ_{k,j}\ \ \text{and}\ \ \gamma_{k,0} e_1 = Q_{k,j}^Tg_k.
\eequationNN
Overall, for generated $(k,j) \in \N{} \times \{0,\dots,n-1\}$, the algorithm considers the trust-region subproblem for a given trust region radius $\delta \in \R{}_{>0}$ defined as
\bequationNN
  \Scal_{k,j}(\delta) := \left[ \min_{t \in \R{j+1}}\ \gamma_{k,0} e_1^Tt + \thalf t^T T_{k,j} t\ \st\ \|t\| \leq \delta \right]
\eequationNN
and/or the regularized subproblem for a given regularization parameter $\lambda \in \R{}_{>0}$ (sufficiently large such that $T_{k,j} + \lambda I \succ 0$) defined as
\bequationNN
  \Rcal_{k,j}(\lambda) := \left[ \min_{t \in \R{j+1}}\ \gamma_{k,0} e_1^Tt + \thalf t^T (T_{k,j} + \lambda I) t \right].
\eequationNN
(We drop the constant objective term $f_k$ in both subproblems since it does not affect the solution sets.)  The key feature of these subproblems is the fact that the matrix~$T_{k,j}$ is tridiagonal, meaning that both can be solved to high accuracy (i.e., exactly for the purposes of our theoretical analysis) in an efficient manner.  In particular, $\Scal_{k,j}(\delta)$ can be solved using the Mor\'e-Sorenson method \cite{MorS83} while $\Rcal_{k,j}(\lambda)$ can be solved by solving the (nonsingular) tridiagonal system $(T_{k,j} + \lambda I) t = -\gamma_{k,0} e_1$.  For future reference, we note that necessary and sufficient conditions for global optimality with respect to $\Scal_{k,j}(\delta)$ are that $(t_{k,j},\lambda_{k,j}) \in \R{j+1} \times \R{}$ is globally optimal if and only if
\bsubequations\label{eq.tltr_opt}
  \begin{align}
    \gamma_{k,0}e_1 + (T_{k,j} + \lambda_{k,j} I)t_{k,j} &= 0, \label{eq.tltr_opt1} \\
    (T_{k,j} + \lambda_{k,j} I) &\succeq 0, \label{eq.tltr_opt2} \\ \text{and}\ \ 
    0 \leq \lambda_{k,j} \perp (\delta - \|t_{k,j}\|) &\geq 0. \label{eq.tltr_opt3}
  \end{align}
\esubequations

Each iteration $k \in \N{}$ of \ITRACE{} begins by computing a solution of $\Scal_{k,j}(\delta_k)$ for some $\delta_k \in \R{}_{>0}$ and sufficiently large $j \in \N{}$.  In particular, \ITRACE{} employs the \emph{truncated Lanczos trust region} algorithm (see \cite[Algorithm~5.2.1]{ConGT00a}) that solves trust-region subproblems over Krylov subspaces of increasing size (i.e., increasing~$j$) until a termination condition holds.  We state our variant of the algorithm (\TLTR) in detail as Algorithm~\ref{alg.tltr}, which generates the aforementioned quantities $\{Q_{k,j}\}_{j\in\N{}}$ and $\{T_{k,j}\}_{j\in\N{}}$ as well as some auxiliary values required for Lanczos.  For our purposes with \ITRACE, the termination conditions that we use are written in the full space as
\bsubequations\label{eq.tltr_termination}
  \begin{align}
    g_k + (H_k + \lambda_{k,j} I)s_{k,j} &= r_{k,j}, \label{eq.tltr_termination_1} \\
    r_{k,j}^Ts_{k,j} &\leq 0, \label{eq.tltr_termination_2} \\
    s_{k,j}^T(H_k + \lambda_{k,j} I) s_{k,j} & \geq 0, \label{eq.tltr_termination_3} \\
    \text{and}\ \ 0 \leq \lambda_{k,j} \perp (\delta_k - \|s_{k,j}\|) &\geq 0 \label{eq.tltr_termination_4}
  \end{align}
\esubequations
along with
\bsubequations\label{eq.tltr_termination_last}
  \begin{align}
    \text{either} && \|r_{k,j}\| &\leq \xi_1 \|s_{k,j}\|^2 \label{eq.tltr_termination_5} \\
    \text{or both} && \|r_{k,j}\| &\leq \xi_2 \min\{1, \|s_{k,j}\|\} \|g_k\| \label{eq.tltr_termination_6} \\
    \text{and} && 1 &\leq \xi_3 \min\{1, \|s_{k,j}\|\} \|T_{k,j} + \lambda_{k,j} I\|, \label{eq.tltr_termination_6_2}
  \end{align}
\esubequations
where $(\xi_1,\xi_2,\xi_3) \in \R{}_{>0} \times (0,1) \times \R{}_{>0}$ are user-prescribed parameters.  The algorithm may also impose tighter residual conditions to achieve a fast rate of local convergence; these are specified along with our analysis in Section~\ref{sec.local_convergence}.

Algorithm~\ref{alg.tltr} generates, for each generated value of $j \in \{0,\dots,n-1\}$, the primal-dual solution $(t_{k,j},\lambda_{k,j})$ of $\Scal_{k,j}(\delta_k)$, the primal solution of which defines the full-space trial step $s_{k,j} \in \R{n}$, which in turn defines the residual vector $r_{k,j} \in \R{n}$ by \eqref{eq.tltr_termination_1}.  If $r_{k,j} = 0$, then \eqref{eq.tltr_termination} is essentially the set of necessary and sufficient conditions for the global minimization of $m_k(s)$ over $s \in \R{n}$ such that $\|s\| \leq \delta_k$, the only difference being the relaxed condition that $s_{k,j}^T(H_k + \lambda_{k,j} I)s_{k,j} \geq 0$ rather than $(H_k + \lambda_{k,j}) \succeq 0$.  We show in our theoretical analysis that, by the construction of Algorithm~\ref{alg.tltr}, the conditions in \eqref{eq.tltr_termination} are satisfied for all generated values of $j$, meaning that the only conditions that need to be checked explicitly are those in \eqref{eq.tltr_termination_last}, the first of which is guaranteed to hold by the iteration when $j$ reaches $n-1$, if not earlier.  Our analysis shows that, in the reduced space, \eqref{eq.tltr_termination_last} is equivalent to
\bsubequations\label{eq.tltr_termination_last_equiv}
  \begin{align}
    \text{either} && \mu_{k,j} &\leq \xi_1 \|t_{k,j}\|^2, \label{eq.tltr_termination_5_equiv} \\
    \text{or both} && \mu_{k,j} &\leq \xi_2 \min\{1, \|t_{k,j}\|\} \gamma_{k,0} \label{eq.tltr_termination_6_equiv} \\
    \text{and} && 1 &\leq \xi_3 \min\{1, \|t_{k,j}\|\} \|T_{k,j} + \lambda_{k,j} I\|, \label{eq.tltr_termination_6_equiv_2}
  \end{align}
\esubequations
which are the conditions that are actually checked in the algorithm.  (For simplicity, \ITRACE{} can always require \eqref{eq.tltr_termination_5_equiv} and simply ignore \eqref{eq.tltr_termination_6_equiv}--\eqref{eq.tltr_termination_6_equiv_2}, but for the sake of generality in our analysis, we show that allowing either condition leads to the same complexity guarantees, which is of interest since, in some situations in practice, the conditions in \eqref{eq.tltr_termination_6_equiv}--\eqref{eq.tltr_termination_6_equiv_2} may be less restrictive than \eqref{eq.tltr_termination_5_equiv}.)

\balgorithm[ht]
  \caption{\TLTR{} (an adaptation of \cite[Algorithm~5.2.1]{ConGT00a})} 
  \label{alg.tltr}
  \balgorithmic[1]
    \renewcommand{\algorithmicrequire}{\textbf{Parameters:}}
    \Require $(\xi_1,\xi_2,\xi_3) \in \R{}_{>0} \times (0,1) \times \R{}_{>0}$ from \ITRACE{} (Algorithm~\ref{alg.itrace})
    \renewcommand{\algorithmicrequire}{\textbf{Input:}}
    \Require $y_{k,0}$ ($= g_k$), $\gamma_{k,0}$ ($= \|g_k\|$), $H_k$, $\delta_k$
    \State $q_{k,-1} \gets 0$, $Q_{k,-1} \gets [\ ]$
    \For{$j = 0,1,\dots$}
      \State $q_{k,j} \gets \frac{1}{\gamma_{k,j}} y_{k,j}$
      \State $Q_{k,j} \gets \bbmatrix Q_{k,j-1} & q_{k,j} \ebmatrix$
      \State $\theta_{k,j} \gets q_{k,j}^T H_k q_{k,j}$ \label{line.tltr.Hq}
      \If{$j=0$}
        \State $T_{k,j} \gets \bbmatrix \theta_{k,0} \ebmatrix$
      \Else
        \State $T_{k,j} \gets \bbmatrix T_{k,j-1} & \begin{matrix} 0 \\ \vdots \\ \gamma_{k,j} \end{matrix} \\ \begin{matrix} 0 & \cdots & \gamma_{k,j} \end{matrix} & \theta_{k,j} \ebmatrix$
      \EndIf
      \State $y_{k,j+1} \gets H_kq_{k,j} - \theta_{k,j}q_{k,j} - \gamma_{k,j}q_{k,j-1}$ \label{line.tltr.y}
      \State $\gamma_{k,j+1} \gets \|y_{k,j+1}\|$
      \State compute $(t_{k,j},\lambda_{k,j})$ by solving $\Scal_{k,j}(\delta_k)$ (which gives $s_{k,j} = Q_{k,j} t_{k,j}$) \label{line.tltr.S}
      \State $\mu_{k,j} \gets \gamma_{k,j+1} |e_{j+1}^Tt_{k,j}|$ ($= \|r_{k,j}\|$ where $r_{k,j} = g_k + (H_k + \lambda_{k,j} I) s_{k,j}$)
      \If{\eqref{eq.tltr_termination_last_equiv} (equivalently, \eqref{eq.tltr_termination_last}) holds} \label{line.tltr.terminate}
        \State \textbf{break}
      \EndIf
    \EndFor
    \State \textbf{return} $(j,t_{k,j},\lambda_{k,j},Q_{k,j},T_{k,j},y_{k,j+1},\gamma_{k,j+1},\mu_{k,j})$
  \ealgorithmic
\ealgorithm

Upon completion of the call to Algorithm~\ref{alg.tltr} in iteration $k \in \N{}$, \ITRACE{} turns to determine whether the current trial solution should be accepted, whether an alternative trial solution should be computed \emph{over the current Krylov subspace} after an expansion and/or contraction(s) of the trust-region radius, or whether the Krylov subspace should be increased in dimension.  This is done by first calling a subroutine that we call ``Find Decrease Step'' (\FDS), stated as Algorithm~\ref{alg.fds} below.  This method checks conditions derived from those in \TRACE{} and ultimately produces a step---potentially after expansion and/or contraction(s) of the trust-region radius---that offers sufficient decrease in the objective and a ratio between the subproblem's dual solution and the norm of the subproblem's primal solution that is sufficiently small, \emph{all while keeping the Krylov subspace fixed}.  The notion of sufficient decrease in the objective, as in \TRACE, uses a function-decrease-to-step-norm ratio of the form
\bequation\label{eq.rho}
  \rho_k(s) := \frac{f_k - f(x_k + s)}{\|s\|^3}.
\eequation
See \cite[Section~2.4]{CurtRobiSama17} for motivation for the use of this ratio for this purpose; in short, using a step acceptance ratio of this form ensures that accepted steps yield a reduction in the objective on the order that is needed to achieve optimal complexity.

\balgorithm[ht]
  \footnotesize
  \caption{Find Decrease Step (\FDS)}
  \label{alg.fds}
  \balgorithmic[1]
    \renewcommand{\algorithmicrequire}{\textbf{Parameters:}}
    \Require $\eta \in (0,1)$, $\underline{\sigma} \in (0, \infty)$, $\overline{\sigma} \in (\underline{\sigma},\infty)$, $\gamma_C \in (0,1)$, and $\gamma_\lambda \in (1,\infty)$ from \ITRACE{} (Algorithm~\ref{alg.itrace})
    \renewcommand{\algorithmicrequire}{\textbf{Input:}}
    \Require $t_{k,j,0}$, $\lambda_{k,j,0}$, $Q_{k,j}$, $T_{k,j}$, $\gamma_{k,0}$, $\gamma_{k,j+1}$, $\mu_{k,j,0}$, $\delta_{k,j,0}$, $\sigma_{k,j,0}$
    \For{$l = 0,1,\dots$}
      \If{$\rho_k(Q_{k,j}t_{k,j,l}) \geq \eta$ and $\frac{\lambda_{k,j,l}}{\|t_{k,j,l}\|} \leq \sigma_{k,j,l}$} \label{line.fds.function_evaluation} \Comment{decrease step found}
        \State \textbf{return} $(t_{k,j,l}, \lambda_{k,j,l}, \mu_{k,j,l}, \delta_{k,j,t}, \sigma_{k,j,l})$
      \ElsIf{$\rho_k(Q_{k,j}t_{k,j,l}) \geq \eta$ and $\frac{\lambda_{k,j,l}}{\|t_{k,j,l}\|} > \sigma_{k,j,l}$} \Comment{expand trust region}
        \State $\delta_{k,j,l+1} \gets \frac{\lambda_{k,j,l}}{\sigma_{k,j,l}}$
        \State compute $(t_{k,j,l+1}, \lambda_{k,j,l+1})$ by solving $\Scal_{k,j}(\delta_{k,j,l+1})$ \label{line.fds.S}
        \State $\sigma_{k,j,l+1} \gets \sigma_{k,j,l}$  \label{line.fds.sigma_expand}
      \Else\ (i.e., $\rho_k(Q_{k,j}t_{k,j,l}) < \eta$) \Comment{contract trust region}
        \If{$\lambda_{k,j,l} < \underline{\sigma} \|t_{k,j,l}\|$} \label{line.fds.contract_not_small}
          \State $\hat{\lambda}_{k,j,l+1} \gets \lambda_{k,j,l} + (\underline{\sigma} \gamma_{k,0})^{1/2}$
          \State compute $\that_{k,j,l+1}$ by solving $\Rcal_{k,j}(\hat{\lambda}_{k,j,l+1})$ \label{line.fds.contract_1}
          \If{$\frac{\hat{\lambda}_{k,j,l+1}}{\|\that_{k,j,l+1}\|} \leq \overline{\sigma}$}
            \State $(t_{k,j,l+1}, \lambda_{k,j,l+1}, \delta_{k,j,l+1}) \gets (\that_{k,j,l+1}, \hat{\lambda}_{k,j,l+1}, \|\that_{k,j,l+1}\|)$ \label{line.fds.contract_return_1}
          \Else (i.e., $\frac{\hat{\lambda}_{k,j,l+1}}{\|\that_{k,j,l+1}\|} > \overline{\sigma}$)
            \State compute $\bar{\lambda}_{k,j,l+1} \in (\lambda_{k,j,l}, \hat{\lambda}_{k,j,l+1})$ so that the solution \label{line.fds.contract_between_1}
            \State \hspace{10pt} $\tbar_{k,j,l+1}$ of $\Rcal_{k,j}(\bar{\lambda}_{k,j,l+1})$ yields $\underline{\sigma} < \frac{\bar{\lambda}_{k,j,l+1}}{\|\tbar_{k,j,l+1}\|} < \overline{\sigma}$ \label{line.fds.contract_between_2}
            \State $(t_{k,j,l+1}, \lambda_{k,j,l+1}, \delta_{k,j,l+1}) \gets (\tbar_{k,j,l+1}, \bar{\lambda}_{k,j,l+1}, \|\tbar_{k,j,l+1}\|)$ \label{line.fds.contract_return_2}
          \EndIf
        \Else\ (i.e., $\lambda_{k,j,l} \geq \underline{\sigma}\|t_{k,j,l}\|$)
          \State $\hat{\lambda}_{k,j,l+1} \gets \gamma_\lambda \lambda_{k,j,l}$ \label{line.fds.contract_not_small_next}
          \State compute $\that_{k,j,l+1}$ by solving $\Rcal_{k,j}(\hat{\lambda}_{k,j,l+1})$ \label{line.fds.solve_linear_system}
          \If{$\|\that_{k,j,l+1}\| \geq \gamma_C \delta_{k,j,l}$}
            \State $(t_{k,j,l+1},\lambda_{k,j,l+1},\delta_{k,j,l+1}) \gets (\that_{k,j,l+1},\hat{\lambda}_{k,j,l+1},\|\that_{k,j,l+1}\|)$ \label{line.fds.contract_return_3}
          \Else\ (i.e., $\|\that_{k,j,l+1}\| < \gamma_C \delta_{k,j,l}$)
            \State $\delta_{k,j,l+1} \gets \gamma_C \delta_{k,j,l}$ \label{line.fds.contract_return_4}
            \State compute $(t_{k,j,l+1}, \lambda_{k,j,l+1})$ by solving $\Scal_{k,j}(\delta_{k,j,l+1})$ \label{line.fds.S2}
          \EndIf
        \EndIf
        \State $\sigma_{k,j,l+1} \gets \max\{\sigma_{k,j,l}, \frac{\lambda_{k,j,l+1}}{\|t_{k,j,l+1}\|}\}$ \label{line.fds.sigma_contract}
      \EndIf
    \State $\mu_{k,j,l+1} \gets \gamma_{k,j+1} |e_{j+1}^Tt_{k,j,l+1}|$
    \EndFor
  \ealgorithmic
\ealgorithm

As previously mentioned, the matrix $T_{k,j}$ in any call to \FDS{} is tridiagonal, meaning that the arising subproblems (i.e., instances of $\Scal_{k,j}$ and/or $\Rcal_{k,j}$) can be solved accurately in an efficient manner.  The only aspect of the algorithm that might raise suspicion is the computation requested in lines~\ref{line.fds.contract_between_1}--\ref{line.fds.contract_between_2}.  However, as for \TRACE{} (see \cite[Appendix]{CurtRobiSama17}), this computation is always well posed (see Lemma~\ref{lem.fds_well_posed} in the next subsection) and is \emph{no more expensive} than solving an instance of $\Scal_{k,j}$ or $\Rcal_{k,j}$.

If the full-space solution corresponding to the output from \FDS{} maintains the desired level of accuracy in the full space (recall~\eqref{eq.tltr_termination_last_equiv}), then---since it has already been shown to yield sufficient decrease and a sufficiently small ratio between the subproblem's dual solution and the norm of the step---\ITRACE{} accepts the step and proceeds to the next iteration.  Otherwise, the dimension of the Krylov subspace is increased---again using the Lanczos process---and \FDS{} is called again to produce a new trial step.  These details can be seen in \ITRACE; see Algorithm~\ref{alg.itrace}.

\balgorithm[ht]
  \caption{\ITRACE}
  \label{alg.itrace}
  \balgorithmic[1]
    \renewcommand{\algorithmicrequire}{\textbf{Parameters:}}
    \Require $(\xi_1,\xi_2,\xi_3) \in \R{}_{>0} \times (0,1) \times \R{}_{>0}$, $\eta \in (0,1)$, $\underline{\sigma} \in (0, \infty)$, $\overline{\sigma} \in (\underline{\sigma},\infty)$, $\gamma_C \in (0,1)$, $\gamma_E \in (1,\infty)$, and $\gamma_\lambda \in (1,\infty)$
    \renewcommand{\algorithmicrequire}{\textbf{Input:}}
    \Require $x_0 \in \R{n}$, $\delta_0 \in (0, \infty)$, $\sigma_0 \in [\underline{\sigma}, \overline{\sigma}]$
    \For{$k = 0,1,\dots$}
      \State $\gamma_{k,0} \gets \|g_k\|$
      \State $(j, t_{k,j}, \lambda_{k,j}, Q_{k,j}, T_{k,j}, y_{k,j+1}, \gamma_{k,j+1}, \mu_{k,j}) \gets \TLTR(g_k, \gamma_{k,0}, H_k, \delta_k)$
      \Loop
        \State $(t_{k,j}, \lambda_{k,j}, \mu_{k,j}, \bar\delta_k, \bar\sigma_k) \gets \FDS(t_{k,j}, \lambda_{k,j}, Q_{k,j}, T_{k,j}, \gamma_{k,0}, \gamma_{k,j+1}, \mu_{k,j}, \delta_k, \sigma_k)$
        \If{\eqref{eq.tltr_termination_last_equiv} (equivalently, \eqref{eq.tltr_termination_last}) holds} \label{line.itrace.residual}
          \State set $s_k \gets Q_{k,j}t_{k,j}$
          \State \textbf{break} \label{line.itrace.break}
        \Else
          \State $j \gets j+1$
          \State $q_{k,j} \gets \frac{1}{\gamma_{k,j}} y_{k,j}$
          \State $Q_{k,j} \gets \bbmatrix Q_{k,j-1} & q_{k,j} \ebmatrix $
          \State $\theta_{k,j} \gets q_{k,j}^TH_kq_{k,j}$ \label{line.itrace.Hq}
          \State $T_{k,j} \gets \bbmatrix T_{k,j-1} & \begin{matrix} 0 \\ \vdots \\ \gamma_{k,j} \end{matrix} \\ \begin{matrix} 0 & \cdots & \gamma_{k,j} \end{matrix} & \theta_{k,j} \ebmatrix$
          \State compute $(t_{k,j},\lambda_{k,j})$ by solving $\Scal_{k,j}(\delta_k)$ (which gives $s_{k,j} = Q_{k,j} t_{k,j}$)
          \State $y_{k,j+1} \gets H_kq_{k,j} - \theta_{k,j}q_{k,j} - \gamma_{k,j}q_{k,j-1}$ \label{line.itrace.y}
          \State $\gamma_{k,j+1} \gets \|y_{k,j+1}\|$
          \State $\mu_{k,j} \gets \gamma_{k,j+1} |e_{j+1}^Tt_{k,j}|$ ($= \|r_{k,j}\|$ where $r_{k,j} = g_k + (H_k + \lambda_{k,j} I) s_{k,j}$)
        \EndIf
      \EndLoop
      \State set $x_{k+1} \gets x_k + s_{k}$
      \State set $\delta_{k+1} \gets \max\{\bar{\delta}_k, \gamma_E \|s_k\|\}$
      \State set $\sigma_{k+1} \gets \bar{\sigma}_k$
    \EndFor
  \ealgorithmic
\ealgorithm

We close this section by noting that \TRACE{} also generates an auxiliary sequence (denoted as $\{\Delta_k\} \subset \R{}_{>0}$ in \cite{CurtRobiSama17}) that influences the step acceptance mechanism, which in the context of \ITRACE{} are the conditions in \FDS{} for determining whether a decrease step has been found or whether the trust-region radius should be expanded or contracted.  The role played by this sequence is to ensure that the algorithm converges from remote starting points even if one were not to assume that the Hessian function is Lipschitz continuous over the path generated by the algorithm iterates.  One could introduce such an auxiliary sequence for \ITRACE{} that would play this same role.  However, as mentioned, for the sake of brevity in this paper, we do not analyze the global convergence properties of the algorithm under this more general setting, and instead have chosen to include upfront (in Assumption~\ref{ass.f}) a Lipschitz continuity assumption for the Hessian function.  In this setting, the auxiliary sequence is not needed to prove the results in this paper, so we have not included it.

\section{Convergence and Complexity Analyses}\label{sec.analysis}

In this section, we prove convergence and worst-case complexity results for \ITRACE{} (Algorithm~\ref{alg.itrace}) under Assumption~\ref{ass.f}.  We also add the following assumption, which is reasonable for the purposes of our analysis since if the algorithm reaches an iteration in which the gradient is zero, then (in a finite number of iterations) it satisfies \eqref{eq.eps_stationary} for any $\epsilon \in \R{}_{>0}$.

\bassumption\label{ass.g_nonzero}
  For all generated $k \in \N{}$, one finds that $g_k \neq 0$.
\eassumption

We begin by proving preliminary results that show that the algorithm is well posed in the sense that it will generate an infinite sequence of iterates.  These results rely heavily on the algorithm's use of the Lanczos process for generating the basis matrices that define the reduced-space subproblems.  We then prove our first main set of results on the algorithm's worst-case complexity properties to approximate first-order stationarity (recall~\eqref{eq.eps_stationary}).  A consequence of these results is that algorithm converges from remote starting points.  Finally, we prove that, as an (inexact) second-order method, \ITRACE{} can achieve a rate of local convergence comparable to \TRACE{}.

In various parts of our analysis, we refer to results in ``standard trust-region theory,'' such as on the relationship between primal and dual trust-region subproblem solutions.  For the sake of brevity, we do not cite particular lemmas for all such results; the reader may refer to textbooks such as \cite{ConGT00,NoceWrig06} for the results that we use.

\subsection{Preliminary results}

In this subsection, we show that \ITRACE{} is well posed in the sense that, for all generated $k \in \N{}$, the \textbf{for} loop in \TLTR{} terminates finitely, each call to \FDS{} terminates finitely, and the inner \textbf{loop} in \ITRACE{} terminates finitely, which together show that \ITRACE{} reaches iteration $k+1$.  Inductively, this means that the algorithm generates iterates \emph{ad infinitum}.

Supposing that \ITRACE{} has reached iteration $k \in \N{}$, our first results in this subsection show that the call to \TLTR{} terminates finitely.  Our presentation of this subroutine is based on the claim that $\mu_{k,j} = \|r_{k,j}\|$ for all generated $j$, where for each such $j$ the vector $r_{k,j}$ is the residual corresponding to $s_{k,j} = Q_{k,j}t_{k,j}$ as defined in~\eqref{eq.tltr_termination_1}.  The following lemma, proved as \cite[Theorem~5.1]{GouLRT99}, formalizes this claim.

\blemma{\cite[Theorem~5.1]{GouLRT99}}\label{lem.tltr_residual}
  For all generated $k \in \N{}$, the call to \TLTR{} yields
  \bequationNN
    r_{k,j} = \left(H_k+\lambda_{k,j}I\right)Q_{k,j}t_{k,j} + g_k = \gamma_{k,j+1} e_{j+1}^Tt_{k,j} q_{k,j+1}\ \text{for all generated}\ j \in \N{},
  \eequationNN
  from which it follows that $\mu_{k,j} \gets \gamma_{k,j+1} |e_{j+1}^Tt_{k,j}| = \|r_{k,j}\|$ for all such $j$.
\elemma

We now show that any call to \TLTR{} by \ITRACE{} terminates finitely.
    
\blemma\label{lem.tltr_finite}
  For all generated $k \in \N{}$, the call to \TLTR{} terminates finitely; more precisely, it terminates in iteration $j$ for some $j \in \{0,\dots,n-1\}$.
\elemma
\bproof
  Consider arbitrary $k \in \N{}$ generated by \ITRACE.  Since it is based on the Lanczos process, it is well known (see, e.g., \cite{ZhaSL17}) that if \TLTR{} continues to iterate, then the dimension of the Krylov subspace $\Kcal_{k,j}$ would increase by 1 whenever $j$ increases by 1 until it reaches some $\jbar \in \{0,\dots,n-1\}$ corresponding to which one finds that
  \bequation\label{eq.gamma_zero}
    \baligned
      \gamma_{k,j} &> 0\ \text{for all $j \in \{0,\dots,\jbar\}$}, \\
      \gamma_{k,\jbar+1} &= 0,\ \ \text{and} \\
      \text{dim}(\Kcal_{k,\jbar}) = \text{dim}(\Kcal_{k,\jbar+1}) &= \jbar + 1.
    \ealigned
  \eequation
  Our goal is to show that \TLTR{} terminates by iteration $\jbar$ at the latest.
  
  Consider arbitrary $j \in \{0,\dots,\jbar\}$.  It follows by Lemma~\ref{lem.tltr_residual} that \eqref{eq.tltr_termination_1} holds.  In addition, by Lemma~\ref{lem.tltr_residual} and the fact that $Q_{k,j}^T q_{k,j+1} = 0$, one finds that
  \bequationNN
    r_{k,j}^Ts_{k,j} = (\gamma_{k,j+1} e_{j+1}^Tt_{k,j} q_{k,j+1})^TQ_{k,j} t_{k,j} = 0,
  \eequationNN
  meaning that \eqref{eq.tltr_termination_2} holds.  Moreover, since $(t_{k,j},\lambda_{k,j})$ is a globally optimal solution of $\Scal_{k,j}(\delta_k)$, it satisfies \eqref{eq.tltr_opt2}, which in turn means that $s_{k,j} \equiv Q_{k,j}t_{k,j}$ and~$\lambda_{k,j}$ yield
  \bequationNN
    \baligned
      s_{k,j}^T(H_k + \lambda_{k,j} I)s_{k,j}
        &= t_{k,j}^T(Q_{k,j}^TH_kQ_{k,j} + \lambda_{k,j} Q_{k,j}^TQ_{k,j})t_{k,j} \\
        &= t_{k,j}^T(T_{k,j} + \lambda_{k,j} I)t_{k,j} \geq 0,
    \ealigned
  \eequationNN
  meaning that \eqref{eq.tltr_termination_3} holds.  Finally, the fact that $(t_{k,j},\lambda_{k,j})$ satisfies \eqref{eq.tltr_opt3} means that $(s_{k,j},\lambda_{k,j})$ (satisfying $\|s_{k,j}\| = \|Q_{k,j}t_{k,j}\| = \|t_{k,j}\|$) satisfies \eqref{eq.tltr_termination_4}.  Overall, it has been shown that \eqref{eq.tltr_termination_1}--\eqref{eq.tltr_termination_4} hold for all $j \in \{0,\dots,\jbar\}$.
  
  To complete the proof, all that remains is to observe that for $j = \jbar$, one finds from Lemma~\ref{lem.tltr_residual} and \eqref{eq.gamma_zero} that $\mu_{k,\jbar} = \|r_{k,\jbar}\| = 0 \leq \xi \|t_{k,j}\|^2$, meaning that \eqref{eq.tltr_termination_5_equiv} holds.  Hence, \TLTR{} would terminate in iteration $\jbar$, if it does not terminate earlier.
\eproof

Our next set of results show that any call to \FDS{} terminates finitely.  First, it is clear that each line of \FDS{} is well posed---since each line involves either a straightforward computation or the computation of a solution of a well-defined subproblem---with the only possible exception being the computation requested in lines~\ref{line.fds.contract_between_1}--\ref{line.fds.contract_between_2}.  The fact that this computation is well posed is proved in the following lemma.

\blemma\label{lem.fds_well_posed}
  For all generated $(k,j,l) \in \N{} \times \N{} \times \N{}$ such that lines~\ref{line.fds.contract_between_1}--\ref{line.fds.contract_between_2} of \FDS{} are reached, the required computation is well posed.
\elemma
\bproof
  Consider arbitrary $(k,j,l) \in \N{} \times \N{} \times \N{}$ such that lines~\ref{line.fds.contract_between_1}--\ref{line.fds.contract_between_2} of \FDS{} are reached.  By construction, one must have $\lambda_{k,j,l}  < \underline{\sigma}\|t_{k,j,l}\|$, $\lambda_{k,j,l} < \hat{\lambda}_{k,j,l+1}$, and $\overline{\sigma} \|\that_{k,j,l+1}\| < \hat{\lambda}_{k,j,l+1}$, where $(t_{k,j,l},\lambda_{k,j,l})$ solves $\Scal_{k,j}(\delta_{k,j,l})$ and $\that_{k,j,l+1}$ solves $\Rcal_{k,j}(\hat\lambda_{k,j,l+1})$.  It follows by standard trust-region theory that the ratio function $\phi : [\lambda_{k,j,l},\infty] \to \Rext{}$ defined by $\phi(\lambda) = \lambda/\|t_{k,j}(\lambda)\|$ with $t_{k,j}(\lambda)$ defined as the solution of $\Rcal_{k,j}(\lambda)$ is monotonically increasing.  Therefore, along with the aforementioned inequalities, it follows that there exists $\bar\lambda_{k,j,l+1}$ such that the solution $\tbar_{k,j,l+1}$ of $\Rcal_{k,j}(\bar\lambda_{k,j,l+1})$ yields $\underline{\sigma} < \frac{\bar{\lambda}_{k,j,l+1}}{\|\tbar_{k,j,l+1}\|} < \overline{\sigma}$, as claimed.
\eproof

Now having shown that each line of \FDS{} is well posed, we proceed to show that the \textbf{for} loop of the subroutine terminates finitely.  We next show that, as is comparable in \TRACE, for all generated $(k,j,l) \in \N{} \times \N{} \times \N{}$ the pair $(t_{k,j,l},\lambda_{k,j,l})$ is a primal-dual globally optimal solution of $\Scal_{k,j}(\delta_{k,j,l})$.  We also show that each such pair corresponds to a pair satisfying \eqref{eq.tltr_termination}, although \eqref{eq.tltr_termination_last_equiv} might not hold.

\blemma\label{lem.fds_opt}
  For all generated $(k,j,l) \in \N{} \times \N{} \times \N{}$, one finds that the pair $(t_{k,j,l}, \lambda_{k,j,l})$ satisfies \eqref{eq.tltr_opt} and the pair $(Q_{k,j}t_{k,j,l},\lambda_{k,j,l})$ $($with $Q_{k,j}t_{k,j,l}$ in place of~$s_{k,j}$ and $\lambda_{k,j,l}$ in place of $\lambda_{k,j}$$)$ satisfies \eqref{eq.tltr_termination}.
\elemma
\bproof
  Consider arbitrary generated $(k,j,l) \in \N{} \times \N{} \times \N{}$.  The desired conclusion about $(t_{k,j,l},\lambda_{k,j,l})$ follows since it is obtained either directly by solving $\Scal_{k,j}(\delta_{k,j,l})$ or by solving $\Rcal_{k,j}(\lambda_{k,j,l})$ to obtain $t_{k,j,l}$, then subsequently setting $\delta_{k,j,l} \gets \|t_{k,j,l}\|$, which again means that $(t_{k,j,l},\lambda_{k,j,l})$ solves $\Scal_{k,j}(\delta_{k,j,l})$.  As for the second desired conclusion, first note from Lemma~\ref{lem.tltr_finite} that it holds for $l=0$.  Then, using the same logic as in the proof of the first desired conclusion, observe for the index $l$ of interest that the pair $(t_{k,j,l},\lambda_{k,j,l})$ satisfies \eqref{eq.tltr_opt} with the same $(Q_{k,j},T_{k,j},\gamma_{k,0},\gamma_{k,j+1})$ as for $l=0$; all that has changed between $l=0$ and the value of $l$ of interest is the trust-region radius.  Hence, the desired conclusion follows using the same logic as in the proof of Lemma~\ref{lem.tltr_finite}, except that \eqref{eq.tltr_termination_5_equiv} (hence, \eqref{eq.tltr_termination_last_equiv}) might no longer be satisfied.
\eproof

Our next lemma shows that each trial step is nonzero, the proof of which is merely an adaptation of \cite[Lemma~3.2]{CurtRobiSama17} to the setting of \ITRACE.

\blemma \label{lem.step_nonzero}
  For all generated $(k,j,l) \in \N{} \times \N{} \times \N{}$, one finds $\|t_{k,j,l}\| > 0$.
\elemma
\bproof
  Consider arbitrary generated $(k,j,l) \in \N{} \times \N{} \times \N{}$.  By Lemma~\ref{lem.fds_opt}, $(t_{k,j,l},\lambda_{k,j,l})$ is a primal-dual globally optimal solution of $\Scal_{k,j}(\delta_{k,j,l})$.  If $T_{k,j} = 0$, then \eqref{eq.tltr_opt1} implies $\gamma_{k,0}e_1 + \lambda_{k,j,l} t_{k,j,l} = 0$, which since $\gamma_{k,0} = \|g_k\| \neq 0$ (under Assumption~\ref{ass.g_nonzero}) means that $\lambda_{k,j,l} \neq 0$ and $t_{k,j,l} \neq 0$, from which the desired conclusion holds.  On the other hand, if $T_{k,j} \neq 0$, then there are two cases.  If $\|t_{k,j,l}\| = \delta_{k,j,l}$, then since $\delta_{k,j,l} > 0$ by construction of the algorithm, the desired conclusion holds; otherwise, $\|t_{k,j,l}\| < \delta_{k,j,l}$, in which case \eqref{eq.tltr_opt1} and \eqref{eq.tltr_opt3} imply $\gamma_{k,0}e_1 + T_{k,j}t_{k,j,l} = 0$, from which it follows (under Assumption~\ref{ass.g_nonzero}) that $\|t_{k,j,l}\| \geq \gamma_{k,0}/\|T_{k,j}\| > 0$.
\eproof

We now show that if a computed dual subproblem solution is sufficiently large relative to the norm of the corresponding primal subproblem solution, then the trust-region constraint must be active and sufficient decrease is offered.
    
\blemma\label{lem.lambda_large}
  For all generated $(k,j,l) \in \N{} \times \N{} \times \N{}$, if $(t_{k,j,l},\lambda_{k,j,l})$ yields
  \bequation\label{eq.lambda_bound}
    \lambda_{k,j,l} \geq (\HLip + 2\eta) \|t_{k,j,l}\|,
  \eequation
  then $\|t_{k,j,l}\| = \delta_{k,j,l}$ and $\rho_k(Q_{k,j}t_{k,j,l}) \geq \eta$.
\elemma
\bproof
  Consider arbitrary generated $(k,j,l) \in \N{} \times \N{} \times \N{}$.  By Lemma~\ref{lem.step_nonzero}, one finds that $\|t_{k,j,l}\| > 0$, implying that $\lambda_{k,j,l} \geq (\HLip + 2\eta) \|t_{k,j,l}\|> 0$.  Hence, by Lemma~\ref{lem.fds_opt} and \eqref{eq.tltr_opt3}, one finds that $\|t_{k,j,l}\| = \delta_{k,j,l}$, which is the first desired conclusion.  Now observe that Assumption~\ref{ass.f} and Taylor's Theorem imply that there exists a point~$\overline{x}_k$ on the line segment $[x_k, x_k + Q_{k,j}t_{k,j,l}]$ such that
  \begin{align*}
    &\ m_k(Q_{k,j}t_{k,j,l}) - f(x_k + Q_{k,j}t_{k,j,l}) \\
    =&\ f_k + g_k^TQ_{k,j}t_{k,j,l} + \thalf t_{k,j,l}^T Q_{k,j}^T H_k Q_{k,j} t_{k,j,l} - f(x_k + Q_{k,j}t_{k,j,l}) \\
    =&\ \thalf t_{k,j,l}^T Q_{k,j}^T (H_k - H(\overline{x}_k)) Q_{k,j}t_{k,j,l} \geq -\thalf \HLip \|t_{k,j,l}\|^3.
  \end{align*}
  On the other hand, since $(Q_{k,j}t_{k,j,l}, \lambda_{k,j,l})$ satisfies \eqref{eq.tltr_termination_1}--\eqref{eq.tltr_termination_3} by Lemma~\ref{lem.fds_opt}, it follows that for some $r_{k,j,l} \in \R{n}$ such that $t_{k,j,l}^TQ_{k,j}^Tr_{k,j,l} \leq 0$ one finds
  \begin{align*}
    &\ f_k - m_k(Q_{k,j}t_{k,j,l}) \\
    =&\ -g_k^TQ_{k,j}t_{k,j,l} - t_{k,j,l}^T Q_{k,j}^T H_k Q_{k,j} t_{k,j,l} + \thalf t_{k,j,l}^T Q_{k,j}^T H_k Q_{k,j} t_{k,j,l} \\
    =&\ -t_{k,j,l}^TQ_{k,j}^T(g_k + H_kQ_{k,j}t_{k,j,l}) + \thalf t_{k,j,l}^T Q_{k,j}^T H_k Q_{k,j} t_{k,j,l} \\
    =&\ -t_{k,j,l}^TQ_{k,j}^T(r_{k,j,l} - \lambda_{k,j,l} Q_{k,j} t_{k,j,l}) + \thalf t_{k,j,l}^T Q_{k,j}^T H_k Q_{k,j} t_{k,j,l} \\
    =&\ -t_{k,j,l}^TQ_{k,j}^Tr_{k,j,l} + \thalf t_{k,j,l}^TQ_{k,j}^T(H_k + \lambda_{k,j,l}I)Q_{k,j}t_{k,j,l} + \thalf \lambda_{k,j,l} \|t_{k,j,l}\|^2 \\
    \geq&\ \thalf \lambda_{k,j,l} \|t_{k,j,l}\|^2.
  \end{align*}
  Therefore, if \eqref{eq.lambda_bound} holds, then one finds that
  \begin{align*}
    \rho_k(Q_{k,j}t_{k,j,l})
    &= \frac{f_k - f(x_k + Q_{k,j}t_{k,j,l})}{\|Q_{k,j}t_{k,j,l}\|^3} \\
    &= \frac{f_k - m_k(Q_{k,j}t_{k,j,l})}{\|t_{k,j,l}\|^3} + \frac{m_k(Q_{k,j}t_{k,j,l}) - f(x_k + Q_{k,j}t_{k,j,l})}{\|t_{k,j,l}\|^3} \\ 
    &\geq \frac{-\thalf \HLip\|t_{k,j,l}\|^3 + \thalf \lambda_{k,j,l} \|t_{k,j,l}\|^2}{\|t_{k,j,l}\|^3} \geq -\thalf \HLip + \thalf (\HLip + 2\eta) = \eta,
  \end{align*}
  as desired.
\eproof

We now partition the set of indices generated within any call to \FDS{} and proceed to show that the number of each type of iteration is finite.  In particular, let us define, for all generated $(k,j) \in \N{} \times \N{}$ such that \FDS{} is called, the index sets
\begin{align*}
  \Acal_{k,j} &:= \left\{l \in \N{} : \text{index $l$ is reached and}\ \rho_k(Q_{k,j}t_{k,j,l}) \geq \eta\ \text{and}\ \frac{\lambda_{k,j,l}}{\|t_{k,j,l}\|} \leq \sigma_{k,j,l} \right\}, \\
  \Ecal_{k,j} &:= \left\{l \in \N{} : \text{index $l$ is reached and}\ \rho_k(Q_{k,j}t_{k,j,l}) \geq \eta\ \text{and}\ \frac{\lambda_{k,j,l}}{\|t_{k,j,l}\|} > \sigma_{k,j,l} \right\}, \\ \text{and}\ \ 
  \Ccal_{k,j} &:= \{l \in \N{} : \text{index $l$ is reached and}\ \rho_k(Q_{k,j}t_{k,j,l}) < \eta\},
\end{align*}
which respectively represent the indices corresponding to \emph{accepted}, \emph{expansion}, and \emph{contraction} steps in the call to \FDS{} corresponding to $(k,j)$.  It follows trivially by construction of \FDS{} that $|\Acal_{k,j}| \leq 1$.  Our next lemma shows that if an expansion or contraction step occurs, then the subsequent step cannot be an expansion step.  This is a critical feature of \TRACE{} as well; see \cite[Lemma~3.7]{CurtRobiSama17}.

\blemma\label{lem.no_expansion}
  For all generated $(k,j,l) \in \N{} \times \N{} \times \N{}$, if $l \in \Ecal_{k,j} \cup \Ccal_{k,j}$, then $(k,j,l+1)$ is generated and $(l+1) \notin \Ecal_{k,j}$.
\elemma
\bproof
  Consider arbitrary generated $(k,j,l) \in \N{} \times \N{} \times \N{}$ such that $l \in \Ecal_{k,j} \cup \Ccal_{k,j}$.  It follows by construction of \FDS{} that iteration $l+1$ will be reached.
  \benumerate
    \item[Case 1:] $\lambda_{k,j,l+1} = 0$.  Since $\sigma_{k,j,0} > 0$ by construction of \ITRACE, it follows by lines \ref{line.fds.sigma_expand} and \ref{line.fds.sigma_contract} in \FDS{} that $\sigma_{k,j,l+1} > 0$.  Hence, one finds that $\sigma_{k,j,l+1} > 0 = \frac{\lambda_{k,j,l+1}}{\|t_{k,j,l}\|}$, from which it follows that $(t+1) \notin \Ecal_{k,j}$.
    \item[Case 2:] $\lambda_{k,j,l+1} > 0$.  By Lemma~\ref{lem.fds_opt} and \eqref{eq.tltr_opt3}, it follows that $\|t_{k,j,l+1}\| = \delta_{k,j,l+1}$.  If $l \in \Ecal_{k,j}$, then one finds that $\|t_{k,j,l+1}\| = \delta_{k,j,l+1} = \frac{\lambda_{k,j,l}}{\sigma_{k,j,l}} > \|t_{k,j,l}\| = \delta_{k,j,l}$, so by standard trust-region theory one finds $\lambda_{k,j,l+1} < \lambda_{k,j,l}$.  Hence,
    \begin{align*}
      \frac{\lambda_{k,j,l+1}}{\|t_{k,j,l+1}\|} \leq \frac{\lambda_{k,j,l}\sigma_{k,j,l}}{\lambda_{k,j,l}} = \sigma_{k,j,l} = \sigma_{k,j,l+1},
    \end{align*}
    from which it follows that $(t+1) \notin \Ecal_{k,j}$.  On the other hand, if $l \in \Ccal_{k,j}$, then line~\ref{line.fds.sigma_contract} ensures that $\sigma_{k,j,l+1} \geq \frac{\lambda_{k,j,l+1}}{\|t_{k,j,l+1}\|}$, so $(l+1) \notin \Ecal_{k,j}$.
  \eenumerate
  The conclusion follows by combining the results of the two cases.
\eproof

An immediate consequence of the previous lemma is that the number of expansion steps in any call to \FDS{} is limited by one.
    
\blemma\label{lem.one_expansion}
  For all generated $(k,j) \in \N{} \times \N{}$ such that \FDS{} is called, $|\Ecal_{k,j}| \leq 1$.
\elemma
\bproof
  Consider arbitrary generated $(k,j) \in \N{} \times \N{}$ such that \FDS{} is called.  If $|\Ecal_{k,j}| = 0$, then there is nothing left to prove.  Otherwise, for some smallest generated $l \in \N{}$ one finds that $l \in \Ecal_{k,j}$.  It then follows by induction that $|\Ecal_{k,j}| = 1$.  After all, since $l \in \Ecal_{k,j}$, Lemma~\ref{lem.no_expansion} implies that $(l+1) \in \Acal_{k,j} \cup \Ccal_{k,j}$.  If $(l+1) \in \Acal_{k,j}$, then \FDS{} terminates and $|\Ecal_{k,j}| = 1$, while if $(l+1) \in \Ccal_{k,j}$, then Lemma~\ref{lem.no_expansion} implies that $(l+2) \in \Acal_{k,j} \cup \Ccal_{k,j}$.  This argument shows inductively that $|\Ecal_{k,j}| = 1$, as claimed.
\eproof

All that remains in order to prove that any call to \FDS{} terminates finitely is to prove that the number of contraction steps is finite.  Toward this end, we now prove that as the result of any contraction step, the trust-region radius is decreased and the dual subproblem solution does not decrease.  The proof of the following lemma is essentially the same as that for \cite[Lemma 3.4]{CurtRobiSama17}, but we provide it for completeness.
    
\blemma\label{lem.contraction_monotonic}
  For all generated $(k,j,l) \in \N{} \times \N{} \times \N{}$, if $l \in \Ccal_{k,j}$, then $(k,j,l+1)$ is generated, $\delta_{k,j,l+1} < \delta_{k,j,l}$, and $\lambda_{k,j,l+1} \geq \lambda_{k,j,l}$.
\elemma
\bproof
  Consider arbitrary generated $(k,j,l) \in \N{} \times \N{} \times \N{}$ such that $l \in \Ccal_{k,j}$.  If line \ref{line.fds.contract_return_1}, \ref{line.fds.contract_return_2}, or \ref{line.fds.contract_return_3} is reached, then $\delta_{k,j,l+1} \gets \|t_{k,j,l+1}\|$ where $t_{k,j,l+1}$ solves $\Rcal_{k,j}(\lambda)$ for $\lambda > \lambda_{k,j,l}$.  Since $\lambda > \lambda_{k,j,l}$, it follows by standard trust-region theory that
  \bequationNN
    \delta_{k,j,l+1} \gets \|t_{k,j,l+1}\| < \|t_{k,j,l}\| \leq \delta_{k,j,l+1}\ \ \text{and}\ \ \lambda_{k,j,l+1} = \lambda > \lambda_{k,j,l}.
  \eequationNN
  The only other possibility is that line \ref{line.fds.contract_return_4} is reached, in which case one finds $\delta_{k,j,l+1} \gets \gamma_C \|t_{k,j,l}\| < \delta_{k,j,l}$, which by standard trust-region theory implies $\lambda_{k,j,l+1} \geq \lambda_{k,j,l}$.
\eproof

We are now prepared to prove that, in any call to \FDS, the number of contraction steps is finite, which along with previous results shows that \FDS{} terminates finitely.

\blemma\label{lem.fds_finite}
  For all generated $(k,j) \in \N{} \times \N{}$ such that \FDS{} is called, the call to \FDS{} terminates finitely.
\elemma
\bproof
  Consider arbitrary generated $(k,j) \in \N{} \times \N{}$ such that \FDS{} is called. As already observed, by construction of \FDS, it follows that $|\Acal_{k,j}| \leq 1$.  Moreover, by Lemma~\ref{lem.one_expansion}, it follows that $|\Ecal_{k,j}| \leq 1$.  Hence, it remains to prove that $|\Ccal_{k,j}| < \infty$.
  
  In order to derive a contradiction, suppose that $|\Ccal_{k,j}| = \infty$, which along with Lemma~\ref{lem.one_expansion} means that $l \in \Ccal_{k,j}$ for all sufficiently large $l \in \N{}$.  Indeed, we may assume, without loss of generality, that $\Ccal_{k,j} = \N{}$.  Our goal now is to show---using the arguments of \cite[Lemma~3.9]{CurtRobiSama17}---that $\{\delta_{k,j,l}\}_{l\in\N{}} \to 0$ and $\{\lambda_{k,j,l}\}_{l\in\N{}} \to \infty$.  By Lemma~\ref{lem.contraction_monotonic} and the fact that $\{\delta_{k,j,l}\}_{l\in\N{}} \subset \R{}_{>0}$ by construction, it follows that $\{\delta_{k,j,l}\}_{l\in\N{}}$ converges.  If line~\ref{line.fds.contract_return_4} is reached infinitely often, then $\{\delta_{k,j,l}\}_{l\in\N{}} \to 0$ and, by standard trust-region theory, $\{\lambda_{k,j,l}\}_{l\in\N{}} \to \infty$, as desired.  Hence, we may assume that line~\ref{line.fds.contract_return_4} is reached only a finite number of times.  Let us now prove that we may also proceed under the assumption that line~\ref{line.fds.contract_return_2} is only reached a finite number of times.  After all, suppose that for some $l \in \N{}$ one finds that line~\ref{line.fds.contract_return_2} is reached, in which case the algorithm sets $(t_{k,j,l+1},\lambda_{k,j,l+1})$ such that during iteration $(l+1) \in \Ccal_{k,j}$ the condition in line~\ref{line.fds.contract_not_small} will test false, meaning that the algorithm will proceed to line~\ref{line.fds.contract_not_small_next} in iteration $l+1$.  Since, by Lemma~\ref{lem.contraction_monotonic}, $\{\delta_{k,j,l}\}_{l\in\N{}}$ is monotonically decreasing and $\{\lambda_{k,j,l}\}_{l\in\N{}}$ is monotonically nondecreasing, it follows that $\{\lambda_{k,j,l}/\|t_{k,j,l}\|\}_{l\in\N{}}$ is monotonically increasing, which means that the condition in line~\ref{line.fds.contract_not_small} will test false in all subsequent iterations, meaning that line~\ref{line.fds.contract_return_2} is only reached a finite number of times, as claimed.  All that remains in order to prove $\{\delta_{k,j,l}\}_{l\in\N{}} \to 0$ and $\{\lambda_{k,j,l}\}_{l\in\N{}} \to \infty$ is to show that these limits hold under the assumption that line~\ref{line.fds.contract_return_1} or line~\ref{line.fds.contract_return_3} is reached for all $l \geq \lbar$ for some $\lbar \in \N{}$.  Under this assumption, one finds that
  \bequationNN
    \lambda_{k,j,l+1} \geq \min\{\lambda_{k,j,l} + (\underline\sigma \gamma_{k,0})^{1/2}, \gamma_\lambda \lambda_k\}\ \ \text{for all}\ \ l \geq \lbar + 1,
  \eequationNN
  which implies that, in fact, $\{\lambda_{k,j,l}\}_{l\in\N{}} \to \infty$.  According to standard trust-region theory, this shows that $\{\delta_{k,j,l}\}_{l\in\N{}} \to 0$, as desired.
  
  Since it has been shown that $\Ccal_{k,j} = \N{}$ implies that one has $\{\delta_{k,j,l}\}_{l\in\N{}} \to 0$ and $\{\lambda_{k,j,l}\}_{l\in\N{}} \to \infty$, one may now conclude from Lemma~\ref{lem.lambda_large} that $l \notin \Ccal_{k,j}$ for some sufficiently large $l \in \N{}$, which is a contradiction to the fact that $\Ccal_{k,j} = \N{}$.
\eproof

We may now prove our concluding result of this subsection.

\blemma\label{lem.infinite_iterates}
  \ITRACE{} generates an infinite sequence of iterates, where for all generated $(k,j) \in \N{} \times \N{}$ one finds that $j \in \{0,\dots,n-1\}$.
\elemma
\bproof
  The result follows by induction.  Supposing that \ITRACE{} reaches iteration $k \in \N{}$, it follows from Lemma~\ref{lem.tltr_finite} that the call to \TLTR{} terminates finitely with $j \leq n-1$ and it follows from Lemma~\ref{lem.fds_finite} that any call to \FDS{} terminates finitely.  Hence, all that remains it to prove that the \textbf{loop} in \ITRACE{} terminates finitely, since this means that \ITRACE{} reaches iteration $k+1$.  This follows using the same argument as in the proof of Lemma~\ref{lem.tltr_finite}, since if $j$ reaches $\jbar \in \{0,\dots,n-1\}$ such that \eqref{eq.gamma_zero} holds, the output from \FDS{} yields $\mu_{k,\jbar} = 0$, in which case the \textbf{loop} will terminate.
\eproof

\subsection{Worst-Case Complexity}\label{sec.complexity}

Our purpose in this subsection is to prove worst-case complexity bounds pertaining to \ITRACE's pursuit of $\epsilon$-stationarity.  In fact, in this subsection we show upper bounds on the total numbers of iterations, function evaluations, derivative evaluations, and Hessian-vector products that \ITRACE{} may perform at iterates at which, for arbitrary $\epsilon \in (0,1)$, the bound \eqref{sec.literature_review} does not hold.  Since this iteration bound holds for arbitrary $\epsilon \in (0,1)$, it follows immediately that \ITRACE{} converges toward first-order stationarity in the limit, i.e., $\{\|g_k\|\} \to 0$.

Our first lemma of this subsection shows that, as in \TRACE, a contraction step causes the ratio between the dual subproblem solution to the norm of the primal subproblem solution to obey certain iteration-dependent and uniform bounds.

\blemma\label{lem.ratio_increase}
  For all generated $(k,j,l) \in \N{} \times \N{} \times \N{}$, if $l \in \Ccal_{k,j}$, then
  \bequationNN
    \underline{\sigma} \leq \frac{\lambda_{k,j,l+1}}{\|t_{k,j,l+1}\|} \leq \max\left\{\overline{\sigma}, \left(\frac{\gamma_{\lambda}}{\gamma_C}\right)\frac{\lambda_{k,j,l}}{\|t_{k,j,l}\|}\right\} \leq \max\left\{\overline{\sigma}, \left(\frac{\gamma_{\lambda}}{\gamma_C}\right)(\HLip + 2\eta)\right\}.
  \eequationNN
  If, in addition, $\lambda_{k,j,l} \geq \underline{\sigma}\|t_{k,j,l}\|$, then
  \bequationNN
    \frac{\lambda_{k,j,l+1}}{\|t_{k,j,l+1}\|} \geq \min \left\{\gamma_\lambda, \frac{1}{\gamma_C}\right\}\frac{\lambda_{k,j,l}}{\|t_{k,j,l}\|}.
  \eequationNN
\elemma
\bproof
  The proof of the first two desired inequalities follows using the same reasoning as in the proof of \cite[Lemma 3.17]{CurtRobiSama17}, the details of which we omit for the sake of brevity.  The next desired inequality follows from Lemma \ref{lem.lambda_large} and the fact that $l \in \Ccal_{k,j}$ only if $\lambda_{k,j,l} < (\HLip + 2\eta)\|t_{k,j,l}\|$.  Finally, under the additional condition that $\lambda_{k,j,l} \geq \underline{\sigma}\|t_{k,j,l}\|$, the final desired conclusion follows using the same reasoning as in the proof of \cite[Lemma 3.23]{CurtRobiSama17}, where again we omit the details for brevity.
\eproof

We now use the previous lemma to prove a critical upper bound.
    
\blemma\label{lem.sigma_max}
  Defining
  \bequationNN
    \sigma_{\max}:= \max\left\{ \sigma_{0},  \overline{\sigma}, \left(\frac{\gamma_{\lambda}}{\gamma_C}\right)(\HLip + 2\eta)\right\} > 0
  \eequationNN
  it follows for all generated $(k,j,l) \in \N{} \times \N{} \times \N{}$ that $\sigma_{k,j,l} \leq \sigma_{\max}$.
\elemma
\bproof
  We prove the result by induction.  As a base case, consider $k=0$ and the corresponding smallest $j \in \N{}$ such that \FDS{} is called.  Given $k=0$ and such $j \in \N{}$, the call to \FDS{} initializes $\sigma_{0,j,0} = \sigma_0 \leq \sigma_{\max}$.  Now suppose that for arbitrary generated $k \in \N{}$ and the corresponding smallest $j \in \N{}$ such that \FDS{} is called one finds for generated $l \in \N{}$ that $\sigma_{k,j,l} \leq \sigma_{\max}$.  If $l \in \Ecal_{k,j}$, then by line~\ref{line.fds.sigma_expand} one finds that $\sigma_{k,j,l+1} = \sigma_{k,j,l} \leq \sigma_{\max}$.  If $l \in \Ccal_{k,j}$, then by Lemma \ref{lem.ratio_increase} and line~\ref{line.fds.sigma_contract} one finds that
  \begin{align*}
    \sigma_{k,j,l+1} \gets \max\left\{\sigma_{k,j,l}, \frac{\lambda_{k,j,l+1}}{\|t_{k,j,l+1}\|}\right\} \leq \max\left\{\sigma_{k,j,l}, \overline{\sigma}, \left(\frac{\gamma_{\lambda}}{\gamma_C}\right)(\HLip + 2\eta)\right\} \leq \sigma_{\max}.
  \end{align*}
  Finally, if $l \in \Acal_{k,j}$, then either $(i)$ \eqref{eq.tltr_termination_last_equiv} is satisfied, \ITRACE{} proceeds to (outer) iteration $k+1$, and for some smallest corresponding $j \in \N{}$ such that \FDS{} is called one finds $\sigma_{k+1,j,0} = \sigma_{k,j,l} \leq \sigma_{\max}$, or $(ii)$ \eqref{eq.tltr_termination_last_equiv} is not satisfied, \ITRACE{} proceeds to (inner) iteration $j+1$, and $\sigma_{k,j+1,0} = \sigma_{k,j,0} \leq \sigma_{\max}$.  Overall, in all cases, the procedures of \ITRACE{} and \FDS{} ensure that the desired conclusion holds.
\eproof

We have proved in Lemma~\ref{lem.infinite_iterates} that \ITRACE{} generates an infinite sequence of iterates, meaning that it generates an infinite sequence of steps $\{s_k\}$.  For our next result, we prove a critical relationship between the norm of each step and the norm of the gradient of the objective function at the subsequent iterate.  Such a relationship is critical for all of the optimal-complexity methods mentioned in Section~\ref{sec.literature_review}.

\blemma\label{lem.step_norm_gradient}
  For all $k \in \N{}$, the step $s_k$ satisfies
  \bequationNN
    \|s_k\| \geq \(\frac{1 - \xi_2}{\thalf \HLip + \sigma_{\max} + \max\{\xi_1,\xi_2 \gLip\}}\)^{1/2} \|g_{k+1}\|^{1/2}.
  \eequationNN
\elemma
\bproof
  Consider arbitrary $k \in \N{}$.  By construction of \FDS{} and \ITRACE{}, one has at line~\ref{line.itrace.break} of \ITRACE{} (with $\lambda_k \gets \lambda_{k,j}$) that $\lambda_k \leq \sigma_{\max} \|s_k\|$ and
  \bequationNN
    \baligned
      \text{either}\ \ \|g_k + (H_k + \lambda_k I)s_k\| &\leq \xi_1 \|s_k\|^2 \\
      \text{or}\ \ \|g_k + (H_k + \lambda_k I)s_k\| &\leq \xi_2 \min\{1,\|s_k\|\} \|g_k\|.
    \ealigned
  \eequationNN
  Under Assumption~\ref{ass.f}, one finds that
  \bequationNN
    \|g_k\| \leq \|g_{k+1}\| + \|g_{k+1} - g_k\| \leq \|g_{k+1}\| + \gLip \|s_k\|;
  \eequationNN
  hence, either $\|g_k + (H_k + \lambda_k I)s_k\| \leq \xi_1 \|s_k\|^2 \leq \xi_1 \|s_k\|^2 + \xi_2\|g_{k+1}\|$ or
  \bequationNN
    \baligned
      \|g_k + (H_k + \lambda_k I)s_k\|
        &\leq \xi_2 \min\{1,\|s_k\|\} (\|g_{k+1}\| + \gLip \|s_k\|) \\
        &\leq \xi_2 (\|g_{k+1}\| + \gLip \|s_k\|^2).
    \ealigned
  \eequationNN
  Overall, since
  \begin{align*}
    \|g_{k+1}\|
    =&\ \|g(x_k + s_k) - (g_k + (H_k + \lambda_k I)s_k) + (g_k + (H_k + \lambda_k I)s_k)\| \\
    \leq&\ \|g(x_k + s_k) - g_k - H_ks_k\| + \lambda_k \|s_k\| + \|g_k + (H_k + \lambda_k I)s_k\|,
  \end{align*}
  it follows from above and by Lemma~\ref{lem.sigma_max} that under Assumption~\ref{ass.f} one has
  \begin{align*}
    &\ (1 - \xi_2) \|g_{k+1}\| \\
    \leq&\ \left\|\int_0^1 (H(x_k + \tau s_k) - H_k)s_k \text{d}\tau \right\| + \sigma_{\max} \|s_k\|^2 + \max\{\xi_1,\xi_2 \gLip\} \|s_k\|^2 \\ 
    \leq&\ \(\int_0^1 \|(H(x_k + \tau s_k) - H_k\| \text{d}\tau \) \|s_k\| + \sigma_{\max} \|s_k\|^2 + \max\{\xi_1,\xi_2 \gLip\} \|s_k\|^2 \\
    \leq&\ \(\int_0^1 \tau \text{d}\tau\) \HLip \|s_k\|^2 + \sigma_{\max}\|s_k\|^2 + \max\{\xi_1,\xi_2 \gLip\} \|s_k\|^2 \\
    \leq&\ \thalf \HLip \|s_k\|^2 + \sigma_{\max} \|s_k\|^2 + \max\{\xi_1,\xi_2 \gLip\}\|s_k\|^2 \\
    =&\ \(\thalf \HLip + \sigma_{\max} + \max\{\xi_1,\xi_2 \gLip\}\) \|s_k\|^2,
  \end{align*}
  which after rearrangement leads to the desired conclusion.
\eproof

It follows from the preceding lemma that the total number of outer iterations that can be performed by \ITRACE{} at iterates at which the norm of the gradient is above $\epsilon \in (0,1)$ is $\Ocal(\epsilon^{-3/2})$, which in turn means that the total number of gradient evaluations at such iterates is also $\Ocal(\epsilon^{-3/2})$.  This is formalized in our first theorem.

\btheorem\label{th.accepted_steps}
  For arbitrary $\epsilon \in (0,1)$, define for \ITRACE{} the index set
  \bequationNN
    \Kcal(\epsilon) := \{k \in \N{} : \|g_k\| > \epsilon\}.
  \eequationNN
  The total number of elements of $\Kcal(\epsilon)$ is at most
  \begin{align*}
    K(\epsilon) := 1 + \left\lfloor \(\frac{(f_0 - f_{\inf})(\thalf \HLip + \sigma_{\max} + \max\{\xi_1,\xi_2 \gLip\})^{3/2}}{\eta (1 - \xi_2)^{3/2}} \) \epsilon^{-3/2} \right\rfloor.
  \end{align*}
  Hence, the total numbers of ``outer'' iterations and gradient evaluations performed at iterates that are not $\epsilon$-stationary are each $\Ocal(\epsilon^{-3/2})$.
\etheorem
\bproof
  By design of \ITRACE{} and Lemma~\ref{lem.step_norm_gradient}, it follows for all $k \in \N{} \setminus \{0\}$ that
  \bequationNN
    f_{k-1} - f_k \geq \eta \|s_{k-1}\|^3 \geq \eta \(\frac{1 - \xi_2}{\thalf \HLip + \sigma_{\max} + \max\{\xi_1,\xi_2 \gLip\}}\)^{3/2} \|g_k\|^{3/2}.
  \eequationNN
  Since $f$ is bounded below under Assumption~\ref{ass.f} and, by construction, \ITRACE{} ensures that $\{f_k\}$ is monotonically nonincreasing, it follows from this string of inequalities that $|\Kcal(\epsilon)| < \infty$.  Therefore, letting $K_\epsilon \in \N{}$ denote the largest index in $\Kcal(\epsilon)$ and summing the prior inequality through iteration $K_\epsilon$ under Assumption~\ref{ass.f} yields
  \begin{align*}
    f_0 - f_{\inf}
      &\geq f_0 - f_{K_\epsilon} \geq \sum_{k=1}^{K_\epsilon} (f_{k-1} - f_k) \\
      &\geq \sum_{k \in \Kcal(\epsilon)} \eta \(\frac{1 - \xi_2}{\thalf \HLip + \sigma_{\max} + \max\{\xi_1,\xi_2 \gLip\}}\)^{3/2} \|g_k\|^{3/2} \\
      &\geq |\Kcal(\epsilon)| \eta \(\frac{1 - \xi_2}{\thalf \HLip + \sigma_{\max} + \max\{\xi_1,\xi_2 \gLip\}}\)^{3/2} \epsilon^{3/2}.
  \end{align*}
  After rearrangement and accounting for iteration $k=0$, the conclusion follows.
\eproof

Our goal now is to account for Hessian-vector products, then function evaluations.  The former occur by line~\ref{line.tltr.Hq} of \TLTR{} and line~\ref{line.itrace.Hq} of \ITRACE{}, and the latter occur by line~\ref{line.fds.function_evaluation} in \FDS{}.  (Hessian-vector products also appear in line~\ref{line.tltr.y} of \TLTR{} and line~\ref{line.itrace.y} of \ITRACE{}, but since these involve the same products as needed in lines~\ref{line.tltr.Hq} and~\ref{line.itrace.Hq}, respectively, one does not need to account for these products as well.  The products can be stored when first computed and reused as needed.)  Our analysis here borrows from the residual analysis from \cite{gould2020error}.  Importantly, in our analysis of \ITRACE{} and its pursuit of (first-order) $\epsilon$-stationarity, we are able to make use of the analysis from~\cite{gould2020error} without having to deal with the so-called \emph{hard case} when solving trust-region subproblems.  This follows from the fact that the worst-case complexity properties for which \ITRACE{} has been designed are of a type described in \cite{gould2020error}, namely, that do not necessitate approximately globally optimal solutions of the arising subproblems.  Indeed, as can be seen in the proof of Lemma~\ref{lem.step_norm_gradient} above, finding subproblem solutions with residuals that are sufficiently small is all that is needed for our purposes.

Following \cite[Section~3.2]{gould2020error}, we note the following.

\blemma\label{lem.T_positive_definite}
  For all generated $(k,j) \in \N{} \times \N{}$, one finds $T_{k,j} + \lambda_{k,j} I \succ 0$.
\elemma
\bproof
  For arbitrary generated $(k,j) \in \N{} \times \N{}$, the conclusion is well known as described in \cite[Section~3.2]{gould2020error}, where it is important to note that, by construction and Lemma~\ref{lem.fds_opt}, the real number $\lambda_{k,j} \in \R{}_{\geq0}$ corresponds to a globally optimal solution of $\Scal_{k,j}(\delta)$ for some $\delta \in \R{}_{\geq0}$ (either $\delta \equiv \delta_k$ or $\delta \equiv \bar\delta_k$ from \FDS).
\eproof

For each generated $(k,j) \in \N{} \times \N{}$, let us write the spectral decomposition
\bequationNN
  T_{k,j} = V_{k,j} \Omega_{k,j} V_{k,j}^T,
\eequationNN
where $V_{k,j} \in \R{(j+1) \times (j+1)}$ is an orthonormal matrix of eigenvectors and $\Omega_{k,j} \in \R{(j+1) \times (j+1)}$ is a diagonal matrix of eigenvalues that are denoted by $\{\omega_{k,j}^{(0)}, \dots, \omega_{k,j}^{(j)})$ and ordered such that $\omega_{k,j}^{(0)} \leq \cdots \leq \omega_{k,j}^{(j)}$.  For all generated $(k,j) \in \N{} \times \N{}$, let us denote the spectral condition number of $T_{k,j} + \lambda_{k,j} I \succ 0$ (recall Lemma~\ref{lem.T_positive_definite}) as
\bequationNN
  \kappa_{k,j} := \frac{\omega_{k,j}^{(j)} + \lambda_{k,j}}{\omega_{k,j}^{(0)} + \lambda_{k,j}} \in \R{}_{>0}.
\eequationNN

Our next result provides an upper bound on the residual defined in \eqref{eq.tltr_termination}.

\blemma\label{lem.residual_convergence}
  For all generated $(k,j) \in \N{} \times \N{}$, one has that
  \bequationNN
    \|r_{k,j}\| \leq \(\frac{2 \|g_k\| H_{\max} \kappa_{k,j}}{\omega_{k,j}^{(j)} + \lambda_{k,j}}\) \(\frac{\sqrt{\kappa_{k,j}} - 1}{\sqrt{\kappa_{k,j}} + 1}\)^j.
  \eequationNN
\elemma
\bproof
  Consider arbitrary generated $(k,j) \in \N{} \times \N{}$.  By construction of \TLTR{} and Lemma~\ref{lem.fds_opt}, it follows that $(t_{k,j},\lambda_{k,j})$ satisfies \eqref{eq.tltr_opt} and $(s_{k,j},\lambda_{k,j})$ with $s_{k,j} = Q_{k,j}t_{k,j}$ satisfies \eqref{eq.tltr_termination}.  Hence, by \cite[Theorem 3.4]{gould2020error}, it follows that
  \bequationNN
    \|r_{k,j}\| \leq \(\frac{2 \|g_k\| \gamma_{k,j+1} \kappa_{k,j}}{\omega_{k,j}^{(j)} + \lambda_{k,j}}\) \(\frac{\sqrt{\kappa_{k,j}} - 1}{\sqrt{\kappa_{k,j}} + 1}\)^j,
  \eequationNN
  and from \cite[Equation~(3.4)]{gould2020error} and Assumption~\ref{ass.f} one finds $\gamma_{k,j+1} \leq \|H_k\| \leq H_{\max}$.  Combining these bounds yields the desired conclusion.
\eproof

We now proceed to prove upper bounds on the total number of inner iterations (over $j \in \N{}$) that are performed during any outer iteration of \ITRACE{} that corresponds to an iterate that is not $\epsilon$-stationary.  For one thing, these bounds serve as upper bounds on the number of Hessian-vector products required during such outer iterations of \ITRACE.  They are also part of upper bounds that we prove for the number of function evaluations during each such outer iteration of \ITRACE.  The first bound that we prove corresponds to the number of iterations that can be performed until \eqref{eq.tltr_termination_5_equiv} holds, whereas the second bound corresponds---assuming \eqref{eq.tltr_termination_6_equiv_2} holds---to the number of iterations that can be performed until \eqref{eq.tltr_termination_6_equiv} holds.  (As has already been seen in the proof of Lemma~\ref{lem.tltr_finite}, the condition in \eqref{eq.tltr_termination_5_equiv} is always satisfiable if enough inner iterations are performed, whereas satisfaction of \eqref{eq.tltr_termination_6_equiv}--\eqref{eq.tltr_termination_6_equiv_2} is not always guaranteed.  That said, the algorithm considers \eqref{eq.tltr_termination_6_equiv}--\eqref{eq.tltr_termination_6_equiv_2} as termination criteria since satisfaction of these inequalities might allow the algorithm to proceed after fewer inner iterations than would be required for \eqref{eq.tltr_termination_5_equiv}.)

To state and prove the aforementioned desired bounds, we define two sets.  Specifically, for arbitrary $(\bar\kappa,\bar\lambda) \in \R{}_{>0} \times \R{}_{>0}$, let us define the sets of index pairs
\bequationNN
  \baligned
    \Ical_1(\bar\kappa,\bar\lambda) := \{(k,j) \in \N{} \times \N{} :&\ \text{$(k,j)$ is generated},\ \kappa_{k,j} \leq \bar\kappa,\ \text{and}\ \omega_{k,j}^{(j)} + \lambda_{k,j} \leq \bar\lambda\} \\ \text{and}\ \ 
    \Ical_2(\bar\kappa) := \{(k,j) \in \N{} \times \N{} :&\ \text{$(k,j)$ is generated and } \kappa_{k,j} \leq \bar\kappa\}.
  \ealigned
\eequationNN
The next lemma shows, at any iterate that is not $\epsilon$-stationary, that if there exists a pair $(\bar\kappa,\bar\lambda)$ such that $(k,j) \in \Ical_1(\bar\kappa,\bar\lambda)$ for sufficiently large $j$, then \TLTR{} and the \textbf{loop} of \ITRACE{} terminate before or at inner iteration number $j$.  It also shows the same conclusion under similar conditions when $(k,j) \in \Ical_2(\bar\kappa)$ and \eqref{eq.tltr_termination_6_equiv_2} holds.  As shown after the lemma, a consequence of this result is that, under nice circumstances including well-conditioning of the (explicitly or implicitly) regularized reduced-space Hessian, the number of iterations performed by \TLTR{} plus the number of iterations of the \textbf{loop} in \ITRACE{} is $\Ocal(\log(\epsilon^{-1}))$.  Otherwise, this sum is at most $\Ocal(n)$.

\blemma\label{lem.product_bound}
  For arbitrary $\epsilon \in (0,1)$ and $k \in \N{}$ such that $\|g_k\| > \epsilon$, consider the following possible scenarios.
  \benumerate
    \item[(i)] There exists $(\bar\kappa,\bar\lambda) \in \R{}_{>0} \times \R{}_{>0}$ such that with
    \bequationNN
      \Jcal_1(\bar\kappa,\bar\lambda) := \min\left\{ n - 1, \left\lceil \log\( \frac{2 H_{\max} \bar\kappa \bar\lambda}{\xi_1 \epsilon} \) \Big/ \log\( \frac{\sqrt{\bar\kappa} + 1}{\sqrt{\bar\kappa} - 1} \) \right\rceil \right\}
    \eequationNN
    one finds that if $(k,j)$ with $j = \Jcal_1(\bar\kappa,\bar\lambda))$ is generated, then $(k,j) \in \Ical_1(\bar\kappa,\bar\lambda)$.
    \item[(ii)] There exists $\bar\kappa \in \R{}_{>0}$ such that with
    \bequationNN
      \Jcal_2(\bar\kappa) := \min\left\{ n - 1, \left\lceil \log\( \frac{2 H_{\max} \bar\kappa \xi_3}{\xi_2 \epsilon} \) \Big/ \log\( \frac{\sqrt{\bar\kappa} + 1}{\sqrt{\bar\kappa} - 1} \) \right\rceil \right\}
    \eequationNN
    one finds that if $(k,j)$ with $j = \Jcal_2(\bar\kappa))$ is generated, then $(k,j) \in \Ical_2(\bar\kappa)$ and $1 \leq \xi_3 \min\{1,\|t_{k,j}\|\} (\omega_{k,j}^{(j)} + \lambda_{k,j})$.
  \eenumerate
  If scenario $(i)$ $($resp.,~$(ii)$$)$ occurs, then \TLTR{} and the \textbf{loop} of \ITRACE{} each terminate before or at inner iteration $j = \Jcal_1(\bar\kappa,\bar\lambda)$ $($resp.,~$j = \Jcal_2(\bar\kappa)$$)$.
\elemma
\bproof
  Consider arbitrary $k \in \N{}$ with $\|g_k\| > \epsilon$.  That each generated $(k,j)$ has $j \leq n-1$ follows from Lemma~\ref{lem.infinite_iterates}.  Hence, all that remains is to prove under the conditions of $(i)$ that each generated $(k,j)$ has $j \leq \Jcal_1(\bar\kappa,\bar\lambda)$, and under the conditions of $(ii)$ that each generated $(k,j)$ has $j \leq \Jcal_2(\bar\kappa)$.
  
  First, suppose the conditions of $(i)$ hold in the sense that either each generated $(k,j)$ has $j < \Jcal_1(\bar\kappa,\bar\lambda)$ or $(k,\Jcal_1(\bar\kappa,\bar\lambda))$ is generated and $(k,\Jcal_1(\bar\kappa,\bar\lambda)) \in \Ical_1(\bar\kappa,\bar\lambda)$.  Observe that $(\sqrt{\kappa} - 1)/(\sqrt{\kappa} + 1)$ is a monotonically increasing function of $\kappa \in \R{}_{>0}$, so by Lemma~\ref{lem.residual_convergence} one finds that, for all generated $(k,j)$, one has
  \bequation\label{eq.residual_bound_kappa}
    \|r_{k,j}\| \leq \(\frac{2 \|g_k\| H_{\max} \bar\kappa}{\omega_{k,j}^{(j)} + \lambda_{k,j}}\) \(\frac{\sqrt{\bar\kappa} - 1}{\sqrt{\bar\kappa} + 1}\)^j.
  \eequation
  Consider the case that $\Jcal_1(\bar\kappa,\bar\lambda) < n-1$ and $j = \Jcal_1(\bar\kappa,\bar\lambda)$, one finds that
  \begin{align}
    && j &\geq \log\( \frac{2 H_{\max} \bar\kappa \bar\lambda}{\xi_1 \epsilon} \) \Big/ \log\( \frac{\sqrt{\bar\kappa} + 1}{\sqrt{\bar\kappa} - 1} \) \nonumber \\
    \iff && j \log\( \frac{\sqrt{\bar\kappa} + 1}{\sqrt{\bar\kappa} - 1} \) &\geq \log\( \frac{2 H_{\max} \bar\kappa \bar\lambda}{\xi_1 \epsilon} \) \nonumber \\
    \iff && j \log\( \frac{\sqrt{\bar\kappa} - 1}{\sqrt{\bar\kappa} + 1} \) &\leq \log\( \frac{\xi_1 \epsilon}{2 H_{\max} \bar\kappa \bar\lambda} \) \nonumber \\
    \iff && \( \frac{\sqrt{\bar\kappa} - 1}{\sqrt{\bar\kappa} + 1} \)^j &\leq \frac{\xi_1 \epsilon}{2 H_{\max} \bar\kappa \bar\lambda}. \label{eq.j_bound}
  \end{align}
  Now observe that, by \eqref{eq.tltr_opt1}, the Cauchy-Schwarz inequality, and the fact that the 2-norm of a real symmetric matrix is its largest eigenvalue, one finds that
  \bequationNN
    \baligned
      \|g_k\|^2 = \|\gamma_{k,0}e_1\|^2 &= \|(T_{k,j} + \lambda_{k,j} I)t_{k,j}\|^2 \\ &\leq \|T_{k,j} + \lambda_{k,j} I\|^2 \|t_{k,j}\|^2 = (\omega_{k,j}^{(j)} + \lambda_{k,j})^2 \|t_{k,j}\|^2.
    \ealigned
  \eequationNN
  Hence, under the conditions of $(i)$, one finds that \eqref{eq.j_bound} implies
  \bequationNN
    \baligned
      \(\frac{2 \|g_k\| H_{\max} \bar\kappa}{\omega_{k,j}^{(j)} + \lambda_{k,j}}\) \( \frac{\sqrt{\bar\kappa} - 1}{\sqrt{\bar\kappa} + 1} \)^j &\leq \(\frac{\|g_k\| }{\omega_{k,j}^{(j)} + \lambda_{k,j}}\) \(\frac{\xi_1 \epsilon}{\bar\lambda}\) \\
      &\leq \xi_1 \(\frac{\|g_k\|^2 }{(\omega_{k,j}^{(j)} + \lambda_{k,j})^2}\) \(\frac{\omega_{k,j}^{(j)} + \lambda_{k,j}}{\bar\lambda}\) \leq \xi_1 \|t_{k,j}\|^2.
    \ealigned
  \eequationNN
  Along with \eqref{eq.residual_bound_kappa}, this bound shows that such $j$ is sufficiently large such that \eqref{eq.tltr_termination_5_equiv} holds.  Therefore, by the construction of \ITRACE{}, the desired conclusion follows.
  
  Now suppose the conditions of $(ii)$ hold in the sense that either each generated $(k,j)$ has $j < \Jcal_2(\bar\kappa)$ or $(k,\Jcal_2(\bar\kappa))$ is generated, $(k,\Jcal_2(\bar\kappa)) \in \Ical_2(\bar\kappa)$, and (since the spectral norm of a symmetric matrix is equal to its largest eigenvalue) with $j = \Jcal_2(\bar\kappa)$ the inequality in \eqref{eq.tltr_termination_6_equiv_2} holds.  The proof in the previous paragraph applies here as well, except with $\Jcal_2(\bar\kappa)$ in place of $\Jcal_1(\bar\kappa,\bar\lambda)$, so that in the present setting \eqref{eq.j_bound} becomes
  \bequation\label{eq.j_bound_2}
    \( \frac{\sqrt{\bar\kappa} - 1}{\sqrt{\bar\kappa} + 1} \)^j \leq \frac{\xi_2 \epsilon}{2 H_{\max} \bar\kappa \xi_3}.
  \eequation
  Under the conditions of $(ii)$, one finds that \eqref{eq.j_bound_2} implies
  \bequationNN
    \(\frac{2 \|g_k\| H_{\max} \bar\kappa}{\omega_{k,j}^{(j)} + \lambda_{k,j}}\) \( \frac{\sqrt{\bar\kappa} - 1}{\sqrt{\bar\kappa} + 1} \)^j \leq \frac{\xi_2 \epsilon \|g_k\|}{\xi_3(\omega_{k,j}^{(j)} + \lambda_{k,j})} \leq \xi_2 \min\{1, \|t_{k,j}\|\} \|g_k\|.
  \eequationNN
  Along with \eqref{eq.residual_bound_kappa} (which applies here as well) and $\|g_k\| = \gamma_{k,0}$, this bound shows that such $j$ is sufficiently large such that \eqref{eq.tltr_termination_6_equiv} and \eqref{eq.tltr_termination_6_equiv_2} hold.  Therefore, by the construction of \ITRACE{}, the desired conclusion follows.
\eproof

We can now prove a worst-case complexity bound for Hessian-vector products.

\btheorem\label{th.products}
  For arbitrary $\epsilon \in (0,1)$, define the index set $\Kcal(\epsilon)$ and positive integer $K(\epsilon)$ as in Theorem~\ref{th.accepted_steps}.  If there exists uniform $(\bar\kappa,\bar\lambda) \in \R{}_{>0} \times \R{}_{>0}$ such that the conditions in $(i)$ and/or $(ii)$ of Lemma~\ref{lem.product_bound} hold for all $k \in \Kcal(\epsilon)$, then the total number of Hessian-vector products performed by \ITRACE{} $($and its subroutines$)$ at iterates that are not $\epsilon$-stationary is at most
  \bequationNN
    K_H(\epsilon) := K(\epsilon) \cdot \min\{\Jcal_1(\bar\kappa,\bar\lambda), \Jcal_2(\bar\kappa)\} = \Ocal(\epsilon^{-3/2} \cdot \min\{n, \log(\epsilon^{-1})\}).
  \eequationNN
  Otherwise, if such $(\bar\kappa,\bar\lambda)$ does not exist, then the number of products is $\Ocal(\epsilon^{-3/2} \cdot n)$.
\etheorem
\bproof
  The result follows by Theorem~\ref{th.accepted_steps} and Lemmas~\ref{lem.infinite_iterates} and \ref{lem.product_bound}.
\eproof

All that remains for our worst-case analysis is to account for function evaluations that occur through line~\ref{line.fds.function_evaluation} in \FDS{}.  Beyond the results that we have proved already, accounting for function evaluations requires proving an upper bound on the number of iterations that can be performed within \FDS.  As is proved in the previous subsection, one finds for all generated $(k,j) \in \N{} \times \N{}$ that $|\Acal_{k,j}| = 1$ and $|\Ecal_{k,j}| \leq 1$ (recall Lemma~\ref{lem.one_expansion}); hence, what is needed for our purposes here is an upper bound on $|\Ccal_{k,j}|$.  A uniform bound over all generated $(k,j) \in \N{} \times \N{}$ is proved in the next lemma.

\blemma\label{lem.contraction_max}
  For all generated $(k,j) \in \N{} \times \N{}$, one finds that
  \begin{align*}
    |\Ccal_{k,j}| \leq 1 + \left\lfloor \frac{\log \(\frac{\sigma_{\max}}{\underline{\sigma}}\)}{\log \(\min \left\{\gamma_\lambda, \frac{1}{\gamma_C}\right\} \) } \right\rfloor =: K_\Ccal.
  \end{align*}
\elemma
\bproof
  Consider arbitrary generated $(k,j) \in \N{} \times \N{}$ such that \FDS{} is called.  If $|\Ccal_{k,j}| = 0$, then the desired conclusion follows trivially.  Hence, we may proceed under the assumption that $|\Ccal_{k,j}| \geq 1$.  It follows by Lemma~\ref{lem.fds_finite} that $|\Ccal_{k,j}| < \infty$.  One may also conclude by Lemma~\ref{lem.no_expansion} that $|\Ccal_{k,j}| \geq 1$ means that $\Ccal_{k,j}$ consists of a set of consecutive positive integers.  Overall, we may proceed knowing that $\Ccal_{k,j} = \{\underline{l},\dots,\overline{l}\}$ for some $(\underline{l},\overline{l}) \in \N{} \times \N{}$.  Since it follows by this definition of $\overline{l}$ and Lemma~\ref{lem.no_expansion} that $(\overline{l}+1) \in \Acal_{k,j}$, it follows with Lemma~\ref{lem.sigma_max} that
  \bequationNN
    \frac{\lambda_{k,j,\overline{l}+1}}{\|t_{k,j,\overline{l}+1}\|} \leq \sigma_{k,j,\overline{l}+1} \leq \sigma_{\max}.
  \eequationNN
  On the other hand, by Lemma~\ref{lem.ratio_increase}, one finds that
  \bequationNN
    \frac{\lambda_{k,j,\overline{l}+1}}{\|t_{k,j,\overline{l}+1}\|} \geq \underline{\sigma} \left(\min \left\{\gamma_\lambda, \frac{1}{\gamma_C}\right\} \right)^{\overline{l} - \underline{l}}.
  \eequationNN
  Combining these upper and lower bounds shows that
  \bequationNN
    (\overline{l} - \underline{l}) \log \(\min \left\{\gamma_\lambda, \frac{1}{\gamma_C}\right\} \) \leq \log\(\frac{\sigma_{\max}}{\underline\sigma}\) \iff \overline{l} - \underline{l} \leq \frac{\log\(\frac{\sigma_{\max}}{\underline\sigma}\)}{\log \(\min \left\{\gamma_\lambda, \frac{1}{\gamma_C}\right\} \)}.
  \eequationNN
  Hence, the desired uniform bound holds since $|\Ccal_{k,j}| = \overline{l} - \underline{l} + 1$.
\eproof

Since, with the previous lemma, there exists a uniform upper bound---independent of the norm of the gradient of the objective---on the number of function evaluations that occur within any inner iteration of \ITRACE, it follows that the worst-case number of function evaluations performed by \ITRACE{} is of the same order as the number of Hessian-vector products.  This is formalized in the following theorem.

\btheorem\label{th.function_evaluations}
  For arbitrary $\epsilon \in (0,1)$, define the index set $\Kcal(\epsilon)$ and positive integer $K(\epsilon)$ as in Theorem~\ref{th.accepted_steps}.  If there exists uniform $(\bar\kappa,\bar\lambda) \in \R{}_{>0} \times \R{}_{>0}$ such that the conditions in $(i)$ and/or $(ii)$ of Lemma~\ref{lem.product_bound} hold for all $k \in \Kcal(\epsilon)$, then the total number of function evaluations performed by \ITRACE{} (and its subroutines) at iterates that are not $\epsilon$-stationary is at most
  \bequationNN
    K_H(\epsilon) := K_\Ccal \cdot K(\epsilon) \cdot \min\{\Jcal_1(\bar\kappa,\bar\lambda), \Jcal_2(\bar\kappa)\} = \Ocal(\epsilon^{-3/2} \cdot \min\{n, \log(\epsilon^{-1})\}).
  \eequationNN
  Otherwise, if such $(\bar\kappa,\bar\lambda)$ does not exist, then the number of evaluations is $\Ocal(\epsilon^{-3/2} \cdot n)$.
\etheorem
\bproof
  The result follows by Theorem~\ref{th.accepted_steps} and Lemmas~\ref{lem.infinite_iterates}, \ref{lem.product_bound}, and \ref{lem.contraction_max}.
\eproof

\subsection{Local Convergence}\label{sec.local_convergence}

\ITRACE{} can attain the same local convergence rate to a strict local minimizer that is attained by \TRACE.  This property of \ITRACE{} follows using well-known results from analyses of inexact Newton methods; nonetheless, it is important to state the results for the sake of completeness.

Our presentation here borrows from that in \cite[Section~3.4]{CurtRobiSama17}.  We consider the local convergence rate attainable by \ITRACE{} under the following assumption.

\bassumption\label{ass.local}
  With respect to an infinite index set $\Scal \subseteq \N{}$, the iterate subsequence $\{x_k\}_{k \in \Scal}$ converges to $x_* \in \R{n}$ at which $H(x_*) \succ 0$.  In addition, there exists a nonempty neighborhood of $x_*$ over which the Hessian function $H$ is locally Lipschitz continuous with Lipschitz constant $\HLoc \in \R{}_{>0}$.
\eassumption

The following lemma captures a property of \TRACE{} inherited by \ITRACE{}.

\blemma\label{lem.local}
  Under Assumption~\ref{ass.local}, the entire sequence $\{x_k\}$ converges to $x_*$.
\elemma
\bproof
  As previously mentioned at the beginning of Section~\ref{sec.complexity}, the analysis in Section~\ref{sec.complexity} shows that $\{\|g_k\|\} \to 0$, which under Assumption~\ref{ass.local} implies $g(x_*) = 0$.  Like in the context of \cite[Lemma~3.30]{CurtRobiSama17}, the remainder of the proof follows similarly to that of \cite[Theorem~6.5.2]{ConGT00a}.
\eproof

Our next lemma is similar to \cite[Lemma~3.31]{CurtRobiSama17} insofar as it shows that, eventually, all computed steps are (potentially inexact) Newton steps that are accepted by the algorithm.  Our proof follows closely that of \cite[Lemma~3.31]{CurtRobiSama17}, but with modifications to account for the potential inexactness of the computed subproblem solutions.

\blemma\label{lem.inactive}
  There exists $k_* \in \N{}$ such that, for all $k \in \N{}$ with $k \geq k_*$, line~\ref{line.itrace.break} of \ITRACE{} is reached with $\lambda_{k,j} = 0$ and $|\Ccal_{k,j}| = |\Ecal_{k,j}| = 0$.
\elemma
\bproof
  By Lemma~\ref{lem.local}, the iterate sequence $\{x_k\}$ converges to $x_*$, at which it follows under Assumption~\ref{ass.local} that $H_* := H(x_*) \succ 0$.  Let the smallest and largest eigenvalues of $H_*$ be denoted by $\omega_{\min}$ and $\omega_{\max}$, respectively.  By continuity of $H$, it follows that the eigenvalues of $H_k$ are contained within the positive interval $[\thalf\omega_{\min},2\omega_{\max}]$ for all sufficiently large $k \in \N{}$.  Consider arbitrary such $k$ and consider arbitrary $j \in \N{}$ such that the index pair $(k,j)$ is generated and \FDS{} is called.  Due to the aforementioned property of the eigenvalues of $H_k$, it follows (see \cite[Section~3.2]{gould2020error}) that the eigenvalues of $T_{k,j}$ are contained in $[\thalf\omega_{\min},2\omega_{\max}]$ as well.  Consider now arbitrary generated $l \in \N{}$.  Either $\|t_{k,j,l}\| = \delta_{k,j,l}$ or $t_{k,j,l} = -T_{k,j}^{-1}(\gamma_{k,0}e_1)$; either way,
  \bequation\label{eq.g_lower_local}
    \|t_{k,j,l}\| \leq \|T_{k,j}^{-1}(\gamma_{k,0}e_1)\| \leq \|T_{k,j}^{-1}\| \|g_k\| \implies \|g_k\| \geq \|t_{k,j,l}\|/\|T_{k,j}^{-1}\|.
  \eequation
  By standard trust-region theory pertaining to Cauchy decrease, it now follows that
  \bequationNN
    \baligned
      f(x_k) - m_k(Q_{k,j}t_{k,j,l})
        &\geq \thalf \|g_k\| \min \left\{ \delta_{k,j,l}, \frac{\|g_k\|}{\|T_{k,j}\|} \right\} \\
        &\geq \thalf \(\frac{\|t_{k,j,l}\|}{\|T_{k,j}^{-1}\|}\) \min \left\{ \|t_{k,j,l}\|, \frac{\|t_{k,j,l}\|}{\|T_{k,j}\|\|T_{k,j}\|^{-1}} \right\} \\
        &\geq \thalf \( \frac{\|t_{k,j,l}\|^2}{\|T_{k,j}\| \|T_{k,j}\|^{-2}} \) \geq \tfrac{1}{16} \omega_{\max}^{-1} \omega_{\min}^2 \|t_{k,j,l}\|^2 =: \eta_* \|t_{k,j,l}\|^2.
    \ealigned
  \eequationNN
  One also finds from \eqref{eq.g_lower_local}, the fact that $\{\|g_k\|\} \to 0$, and the aforementioned properties of the eigenvalues of $T_{k,j}$ that for any $\varepsilon \in (0,1)$ there exists sufficiently large $k_\varepsilon \in \N{}$ such that $\|t_{k,j,l}\| \leq \varepsilon$ for all generated $(k,j,l)$ with $k \geq k_\varepsilon$.  Combining these facts shows, using a similar argument as in the proof of Lemma~\ref{lem.lambda_large}, that for sufficiently large $k \in \N{}$ one finds for any generated $(k,j,l)$ that
  \bequationNN
    \baligned
      f_k - f(x_k + Q_{k,j}t_{k,j,l})
        &\geq f_k - m_k(Q_{k,j}t_{k,j,l}) + m_k(Q_{k,j}t_{k,j,l}) - f(x_k + Q_{k,j}t_{k,j,l}) \\
        &\geq \eta_* \|t_{k,j,l}\|^2 - \thalf \HLoc \|t_{k,j,l}\|^3 \geq \eta \|Q_{k,j}t_{k,j,l}\|^3.
    \ealigned
  \eequationNN
  It follows from this fact that, for any such generated $(k,j,l)$, one has $l \in \Acal_{k,j} \cup \Ecal_{k,j}$.
  
  By the results of the previous paragraph, there exists $\delta_{\min} \in \R{}_{>0}$ such that $\delta_{k,j,l} \geq \delta_{\min}$ for all generated $(k,j,l)$ with sufficiently large $k$.  In addition, continuity of $g$ and the aforementioned properties of the eigenvalues of $T_{k,j}$ imply that the trial step $t_{k,j,l}$ lies in the interior of the trust region for all generated $(k,j,l)$ with sufficiently large $k$.  Since this means that $\lambda_{k,j,l} = 0$ for all such generated $(k,j,l)$, it follows that, in fact, for all generated $(k,j,l)$ for sufficiently large $k$ one has $l \in \Acal_{k,j}$.
\eproof

We now use standard theory of inexact Newton methods to show that \ITRACE{} can, e.g., attain the same rate of local convergence as \TRACE{} (see \cite[Theorem~3.32]{CurtRobiSama17}).

\btheorem
  If, in addition to \eqref{eq.tltr_termination_last_equiv}, the \textbf{if} condition in line \ref{line.itrace.residual} of \ITRACE{} requires $\mu_{k,j} = o(\|g_k\|)$, then $\{x_k\} \to x_*$ Q-superlinearly.  In particular, if the condition requires $\mu_{k,j} = \Ocal(\|g_k\|^2)$, then $\{x_k\} \to x_*$ Q-quadratically.
\etheorem
\bproof
  With Lemma~\ref{lem.local} and \ref{lem.inactive}, the conclusion follows using standard theory of inexact Newton methods; see \cite[Theorem~3.3]{DemES82}.
\eproof

\section{Numerical Results}\label{sec.numerical}

In this section, we provide the results of numerical experiments of a prototype implementation of \ITRACE{}, as well as implementations of \TRACE{} \cite{CurtRobiSama17} and \ARC{} \cite{CartGoulToin11,CartGoulToin11b} for the sake of comparison.  The purposes of presenting these experimental results are twofold.  First, we show that, by allowing inexact subproblem solutions, \ITRACE{} offers computational flexibility beyond that offered by \TRACE{}.  Second, we show that, in terms of key performance measures, \ITRACE{} performs at least as well as \ARC{}, which is a state-of-the-art second-order method that offers optimal complexity to $\epsilon$-stationarity.  For these experiments, all of the algorithms were implemented in a single software package in Matlab.  All experiments were run using the \href{https://coral.ise.lehigh.edu/wiki/doku.php/info:coral}{\texttt{polyps} cluster at Lehigh's COR@L Laboratory}.  Each job was run with a wall-clock-time limit of 90 minutes and a memory limit of 8GB.

\subsection{Implementation details}

The implementations of \ITRACE{} and \TRACE{} share many commonalities.  For a fair comparison, the implementations both involve the auxiliary sequence $\{\Delta_k\}$, the values of which are set and used as in \cite[Algorithm~1]{CurtRobiSama17}.  As explained in the last paragraph of Section~\ref{sec.algorithm}, the theoretical guarantees that have been proved in this paper are maintained with the inclusion of this auxiliary sequence, and in fact allow one to prove guarantees under weaker assumptions.  For our experiments, the common parameters for \ITRACE{} and \TRACE{} were set as $\eta = 10^{-4}$, $\underline{\sigma} = 0.01$, $\overline{\sigma} = 100$, $\gamma_C = 0.5$, $\gamma_E = 1.1$, $\gamma_\lambda = 2$, $\delta_0 = 1$, $\sigma_0 = 1$, and $\Delta_0 = 100$.  Specifically for \ITRACE{}, we ran experiments for $(\xi_1,\xi_2) \in \{(0.1,0.01),(1,0.1),(9,0.9)\}$ and $\xi_3 = 10^6$.  For the implementation of \TRACE{}, all trust-region subproblems are solved using an implementation of the Mor\'e-Sorensen approach \cite{MorS83}.  For the implementation of \ITRACE{}, the $\Rcal_{k,j}$ subproblems are solved by solving a tridiagonal systems, the $\Scal_{k,j}(\delta)$ subproblems are solved using the aforementioned implementation of the Mor\'e-Sorensen approach, and the subproblem in line~\ref{line.fds.contract_between_2} is solved using an implementation of \cite[Algorithm 6.1]{CartGoulToin11}, where, as described in \cite{CurtRobiSama17}, the algorithm is terminated as soon as the ratio $\frac{\bar\lambda_{k,j,l+1}}{\|\tbar_{k,j,l+1}\|}$ lies in the interval $(\underline{\sigma},\overline{\sigma})$.

For the implementation of \ARC, the parameters were set as $\eta_1=10^{-4}$, $\eta_2=0.9$, and $\sigma_0=1$.  In \ARC, $\{\sigma_k\}$ is the sequence of cubic regularization values that is updated dynamically by the algorithm.  In our implementation, this sequence is updated as for the experiments in \cite{CartGoulToin11}, namely, $\sigma_{k+1} \gets \max\{\min\{\sigma_k, \|g_k\|, 10^{-16}\}\}$ if $k$ is a very successful iteration, $\sigma_{k+1} \gets \sigma_k$ if $k$ is successful (but not very successful) iteration, and $\sigma_{k+1} \gets 2\sigma_k$ if $k$ is an unsuccessful iteration.  Like for \ITRACE, the subproblems are solved using an iterative method that employs the Lanczos approach, where for a termination condition our implementation employs \texttt{TC.s} stated as \cite[(3.28)]{CartGoulToin11}, which involves the user-defined parameter $\kappa_\theta$.  Note that \texttt{TC.s} is the same as \eqref{eq.tltr_termination_6_equiv} with $\kappa_\theta \equiv \xi_2$.  Comparable to \ITRACE, we ran experiments with $\kappa_\theta \in \{0.01,0.1,0.9\}$.

All implemented algorithms respect the same termination condition, namely,
\bequation\label{eq.tc_experiment}
  \|g_k\| \leq 10^{-5} \max\{1, \|g_0\|\}.
\eequation

\subsection{Computational flexibility offered by inexactness}\label{sec.versus_trace}

\begin{figure}[ht]
    \centering
    \includegraphics[width=5.2cm]{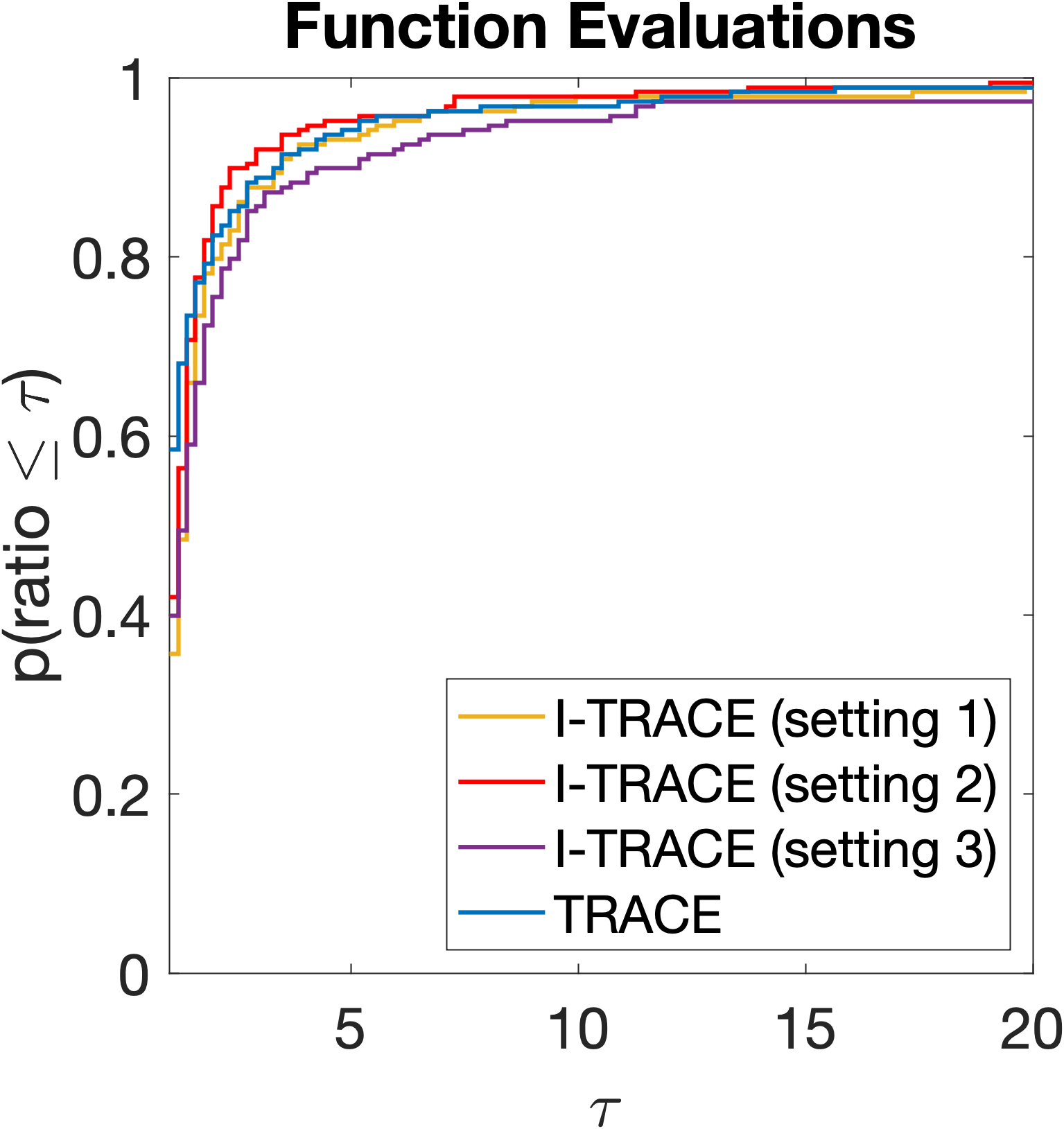}
    \includegraphics[width=5.2cm]{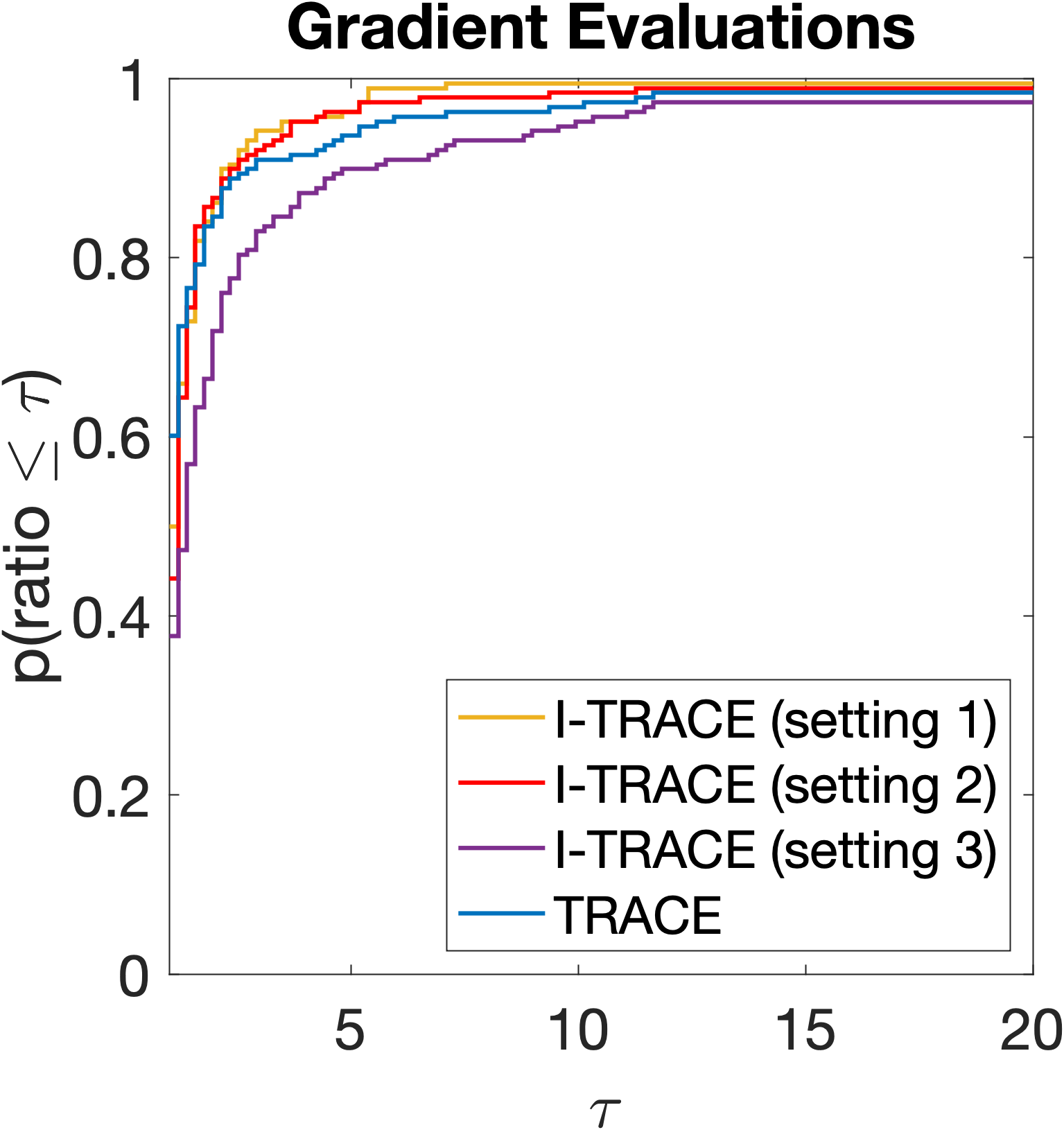}
    \includegraphics[width=5.2cm]{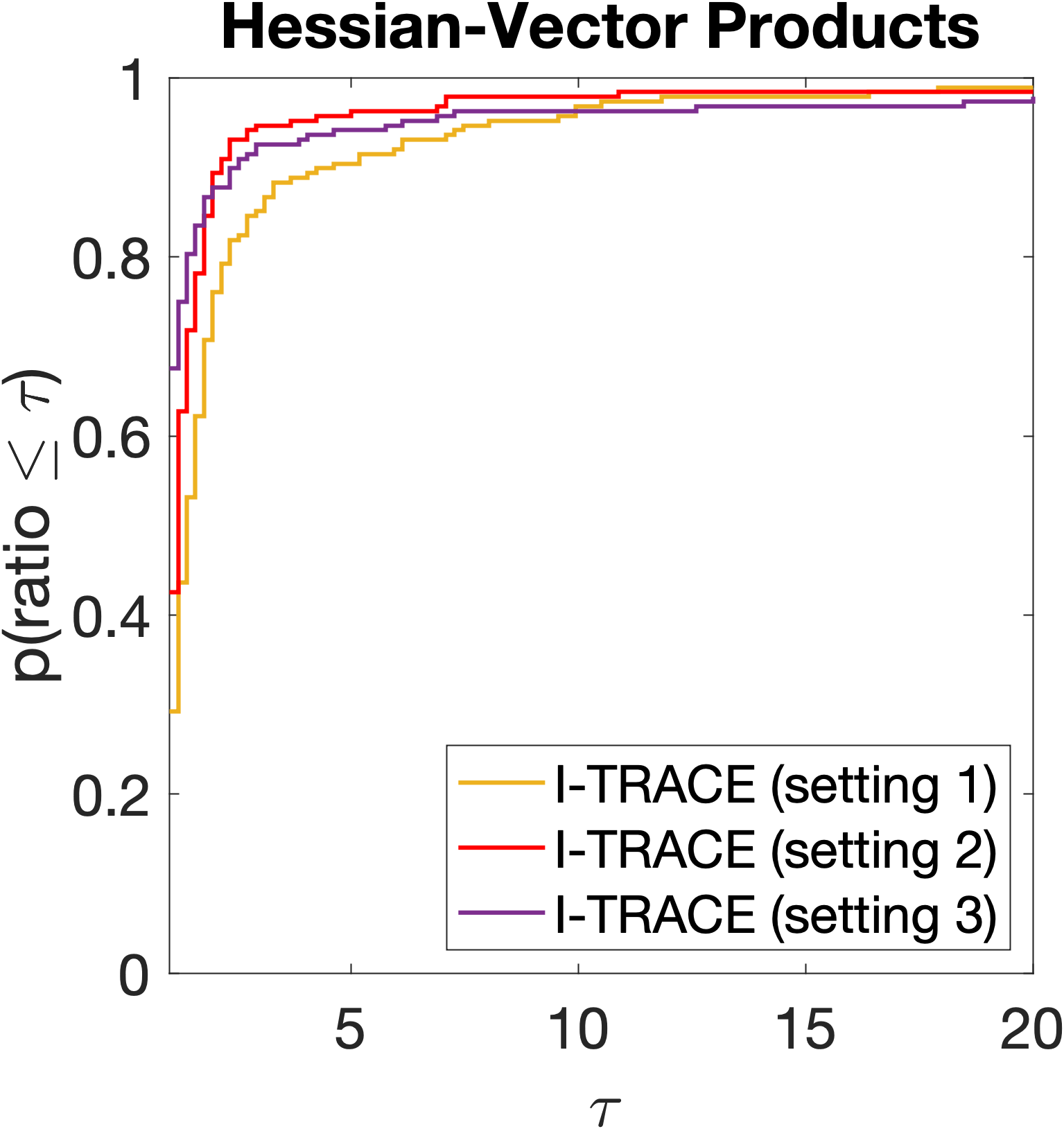}
    \caption{Performance profiles for function evaluations, gradient evaluations, and Hessian-vector products for \ITRACE{} $($with three parameter settings$)$ and \TRACE{} when solving 188 problems from the CUTEst collection.  $($\TRACE{} is not included in the profile for Hessian-vector products since it does not employ Krylov subspace techniques for solving the arising subproblems.$)$} \label{fig.trace_pp}
\end{figure}

Our first set of experiments demonstrates the computational flexibility that \ITRACE{} allows over \TRACE{} due to the fact that \ITRACE{} can employ inexact subproblem solutions.  For this experiment, we ran \ITRACE{} with all parameter settings (see the choices of $(\xi_1,\xi_2)$ in the previous subsection, respectively referred to as ``setting~1,'' ``setting~2,'' and ``setting~3'') and \TRACE{} to solve all of the unconstrained instances in the CUTEst~\cite{GouOT13} collection (with their original parameter settings).  This originally includes 238 problems.  Defining success as encountering an iterate satisfying \eqref{eq.tc_experiment}, \ITRACE{} with setting~1 successfully solved 214 problems, \ITRACE{} with setting~2 successfully solved 218 problems, \ITRACE{} with setting~3 successfully solved 219 problems, and \TRACE{} successfully solved 188 problems (due to hitting the time or memory limit much more often than \ITRACE).  To demonstrate relative performance when solving all problems for which all algorithms/settings were successful (a set of 188 problems), we provide in Figure~\ref{fig.trace_pp} a set of Dolan-Mor\'e performance profiles \cite{DolM02} for function evaluations, gradient evaluations, and Hessian-vector products, respectively.  (We limit the horizontal axis to $\tau=20$ so the differences between the graphs can be seen more clearly.)

The profiles in Figure~\ref{fig.trace_pp} show that, despite allowing inexact subproblem solutions, \ITRACE{} performs comparably to \TRACE{} in terms of function and gradient evaluations, which also means that the algorithms/settings perform comparably in terms of iterations required.  In terms of Hessian-vector products, \ITRACE{} with setting~1 falls a bit behind the other settings, which we contend is due to the algorithm requiring more accurate subproblem solutions in each iteration.  That said, \ITRACE{} with setting~1 performs better in terms of gradient evaluations.  These results demonstrate, as mentioned in Section~\ref{sec.introduction}, that \ITRACE{} offers flexibility between derivative evaluations and Hessian-vector products.  A user can choose the parameters that are preferable depending on the relative costs of these operations for a given problem.


\subsection{Comparison with a state-of-the-art optimal-complexity algorithm}

In this section, we compare the performances of \ITRACE{} and \ARC{}.  First, we mention that \ARC{} with setting~1 successfully solved 211 problems, \ARC{} with setting~2 successfully solved 214 problems, and \ARC{} with setting~3 successfully solved 216 problems; these levels of success were comparable to those for \ITRACE{} (stated in Section~\ref{sec.versus_trace}).

\begin{figure}[ht]
    \centering
    \includegraphics[width=5.2cm]{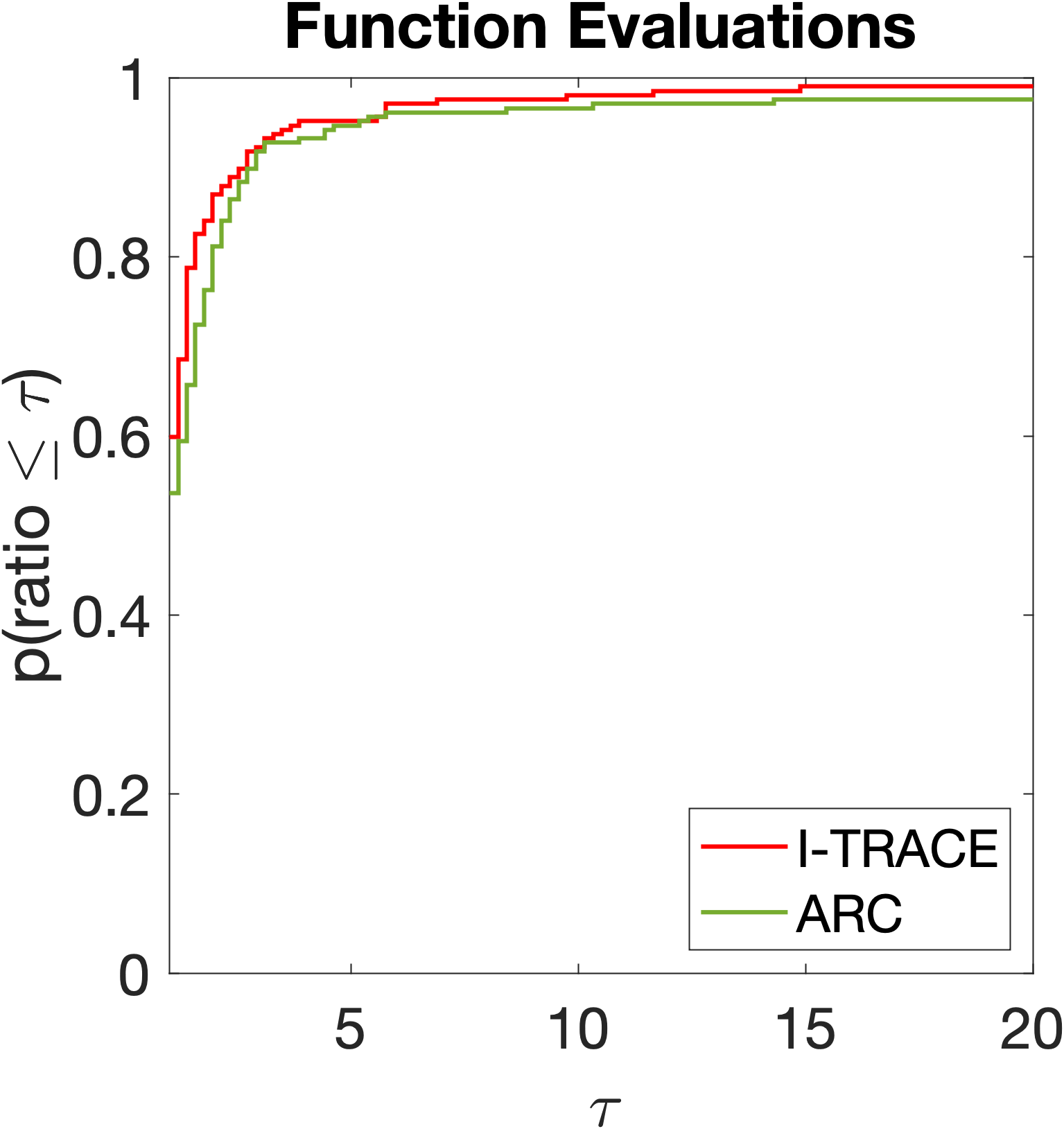}
    \includegraphics[width=5.2cm]{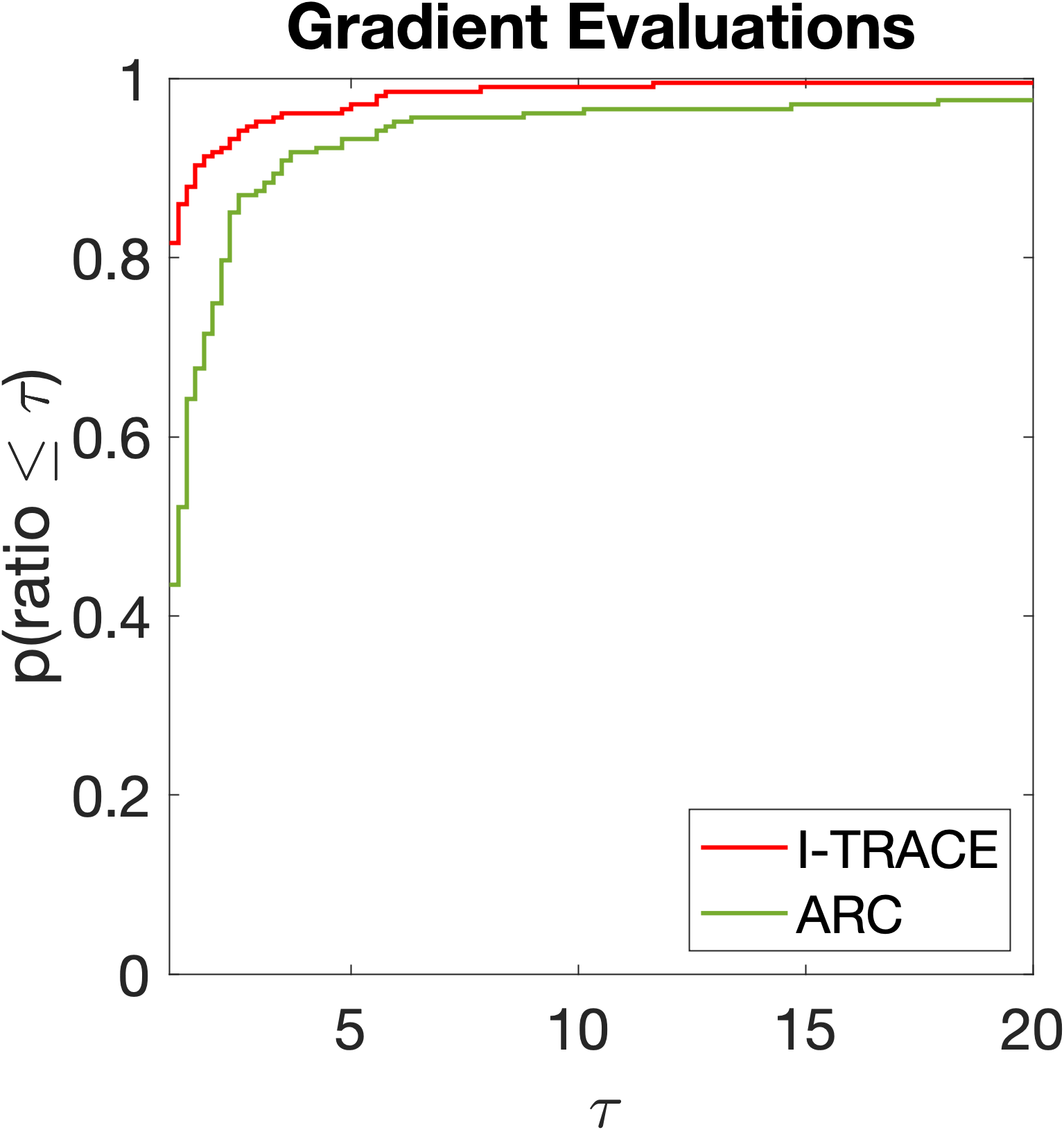}
    \includegraphics[width=5.2cm]{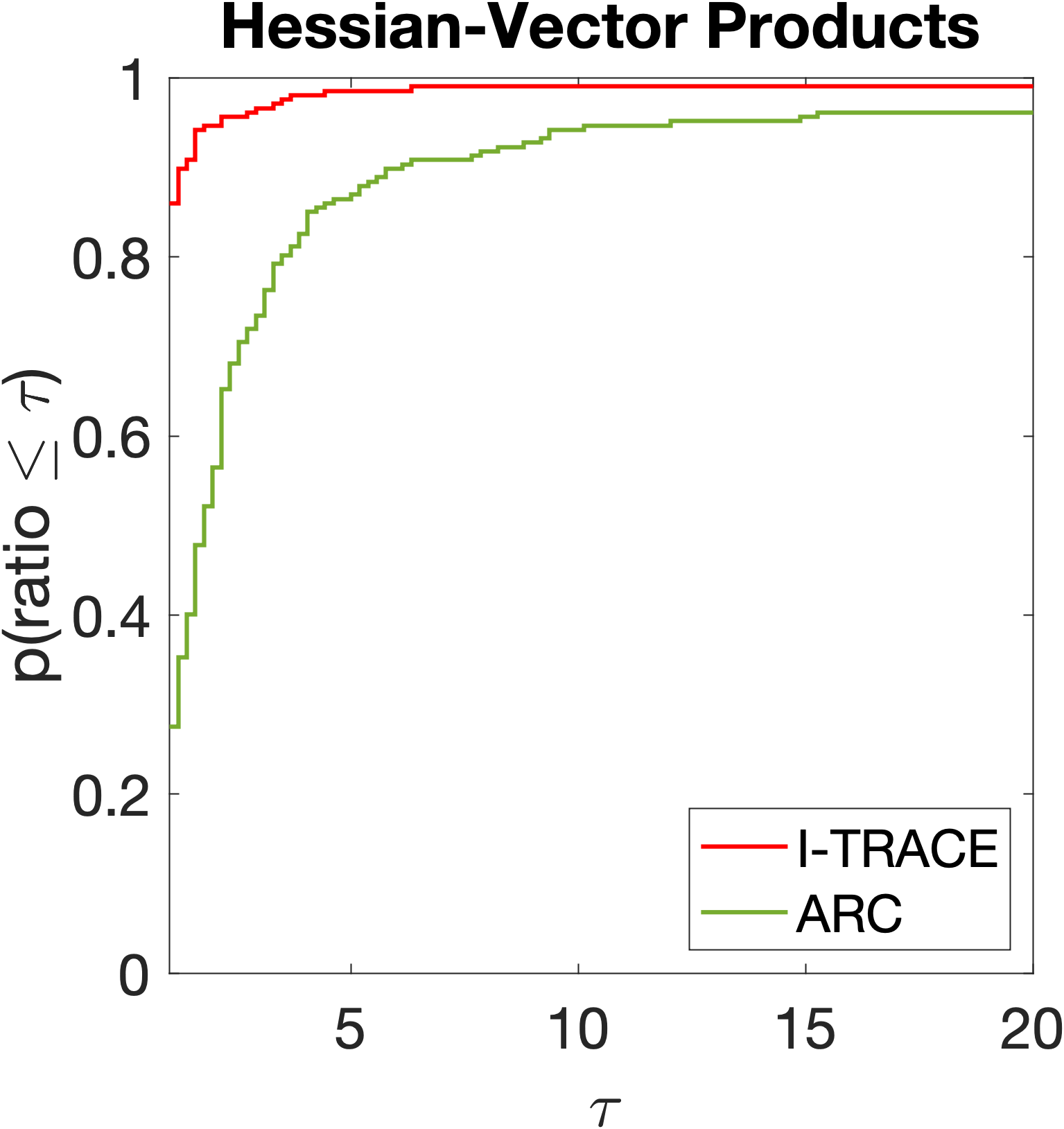}
    \caption{Performance profiles for function evaluations, gradient evaluations, and Hessian-vector products for \ITRACE{} and \ARC{} $($both with setting~1$)$ when solving 207 CUTEst problems.} \label{fig.other_pp1}
\end{figure}

\begin{figure}[ht]
    \centering
    \includegraphics[width=5.2cm]{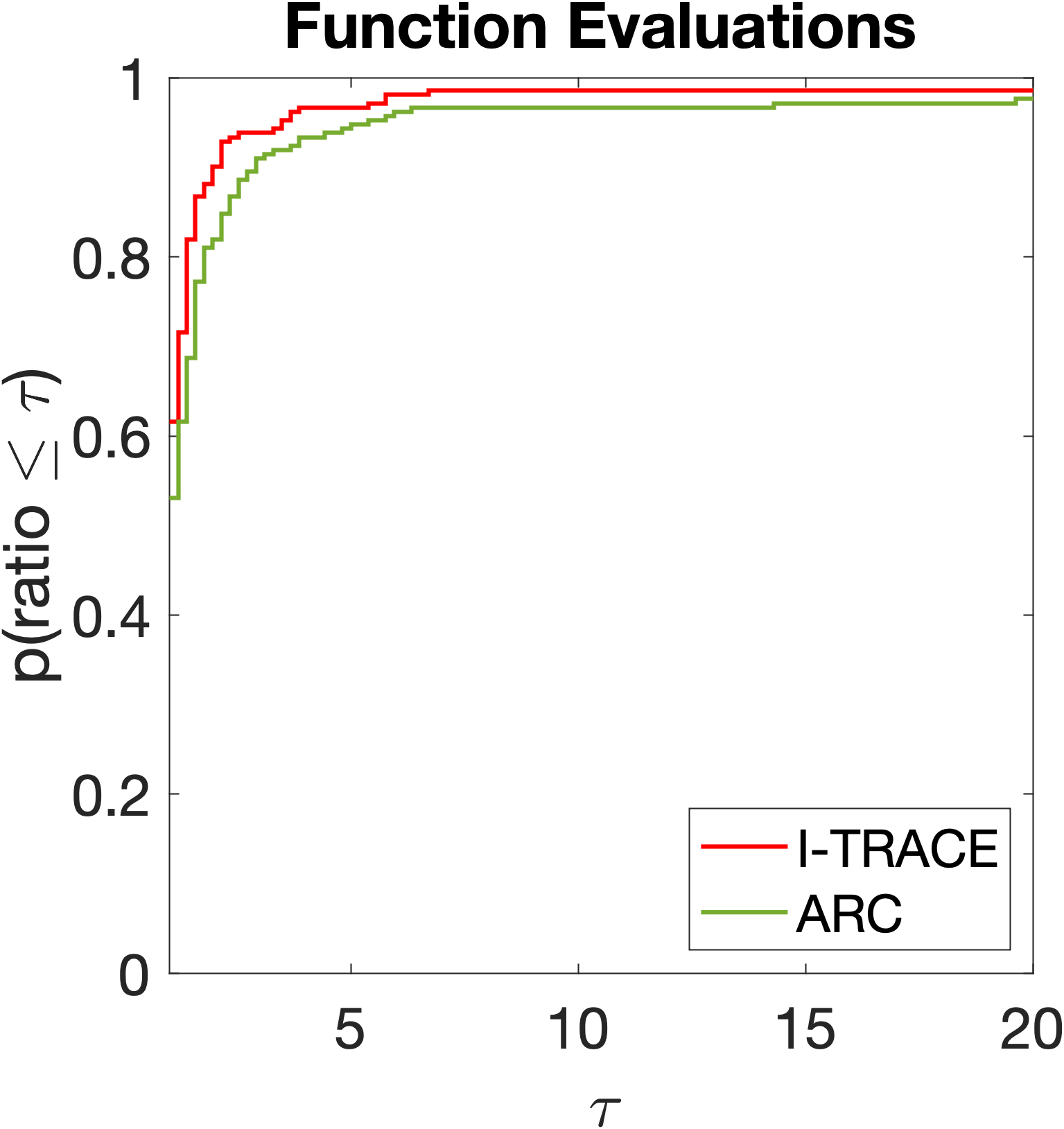}
    \includegraphics[width=5.2cm]{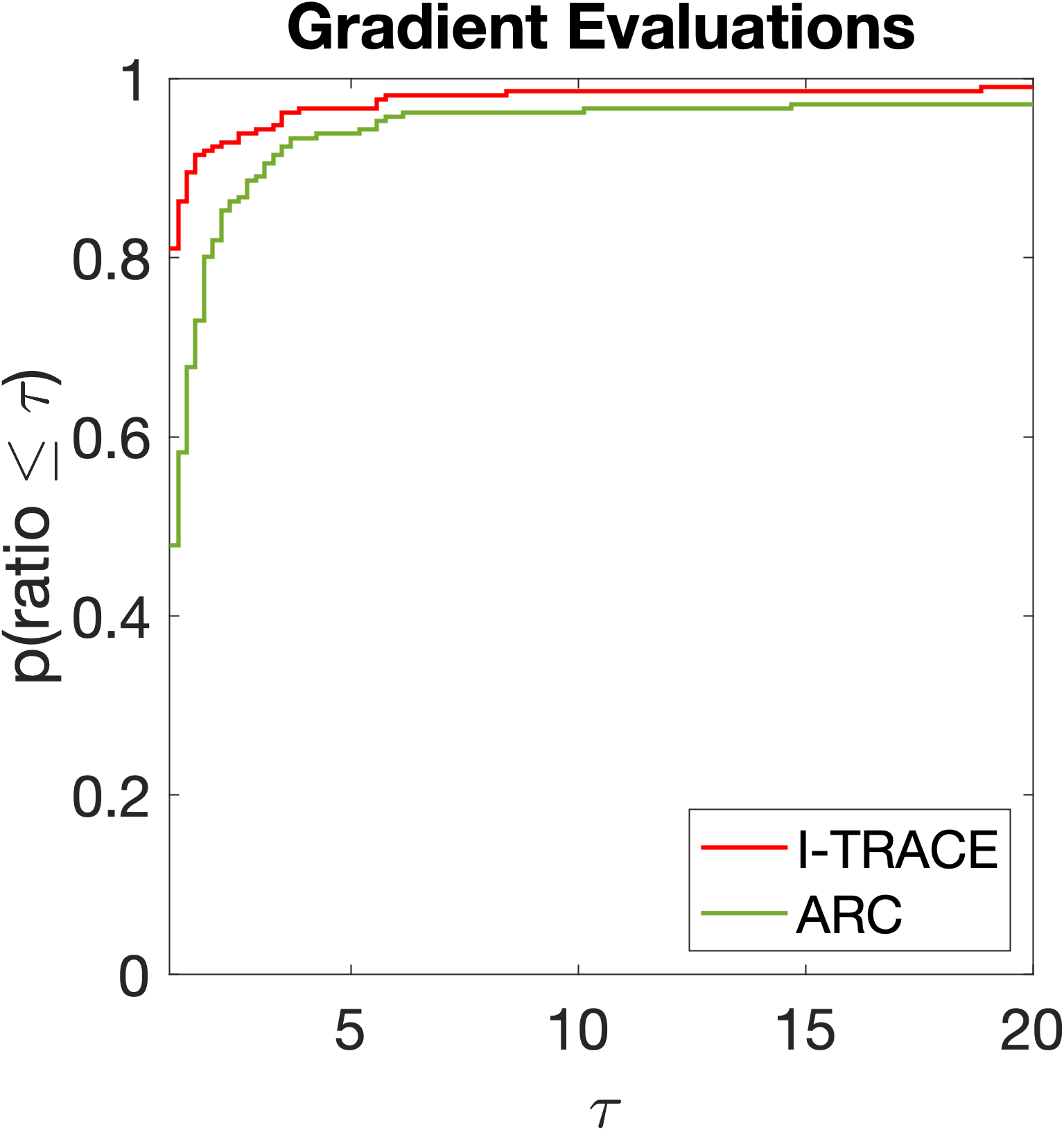}
    \includegraphics[width=5.2cm]{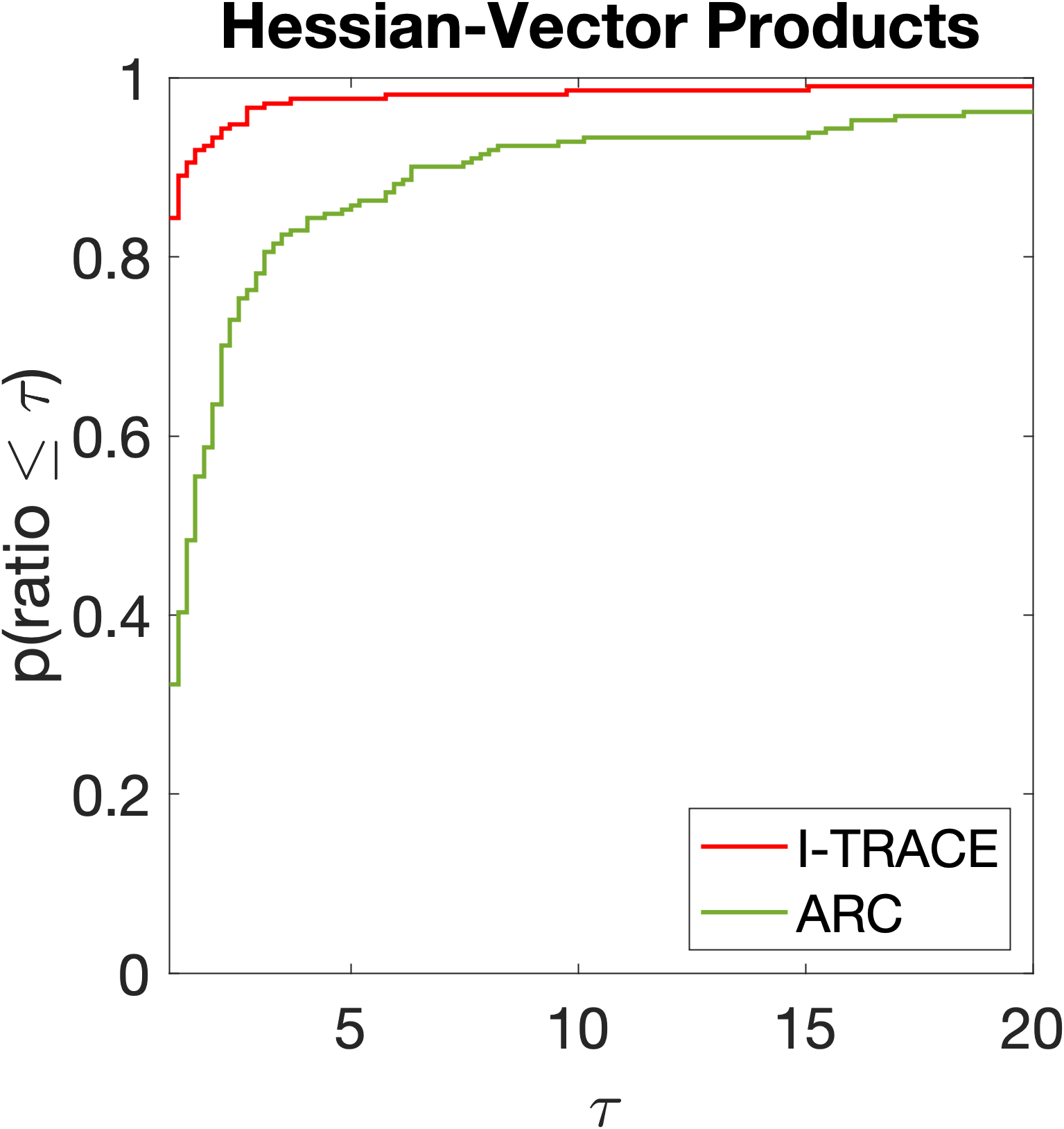}
    \caption{Performance profiles for function evaluations, gradient evaluations, and Hessian-vector products for \ITRACE{} and \ARC{} $($both with setting~2$)$ when solving 211 CUTEst problems.} \label{fig.other_pp2}
\end{figure}

\begin{figure}[ht]
    \centering
    \includegraphics[width=5.2cm]{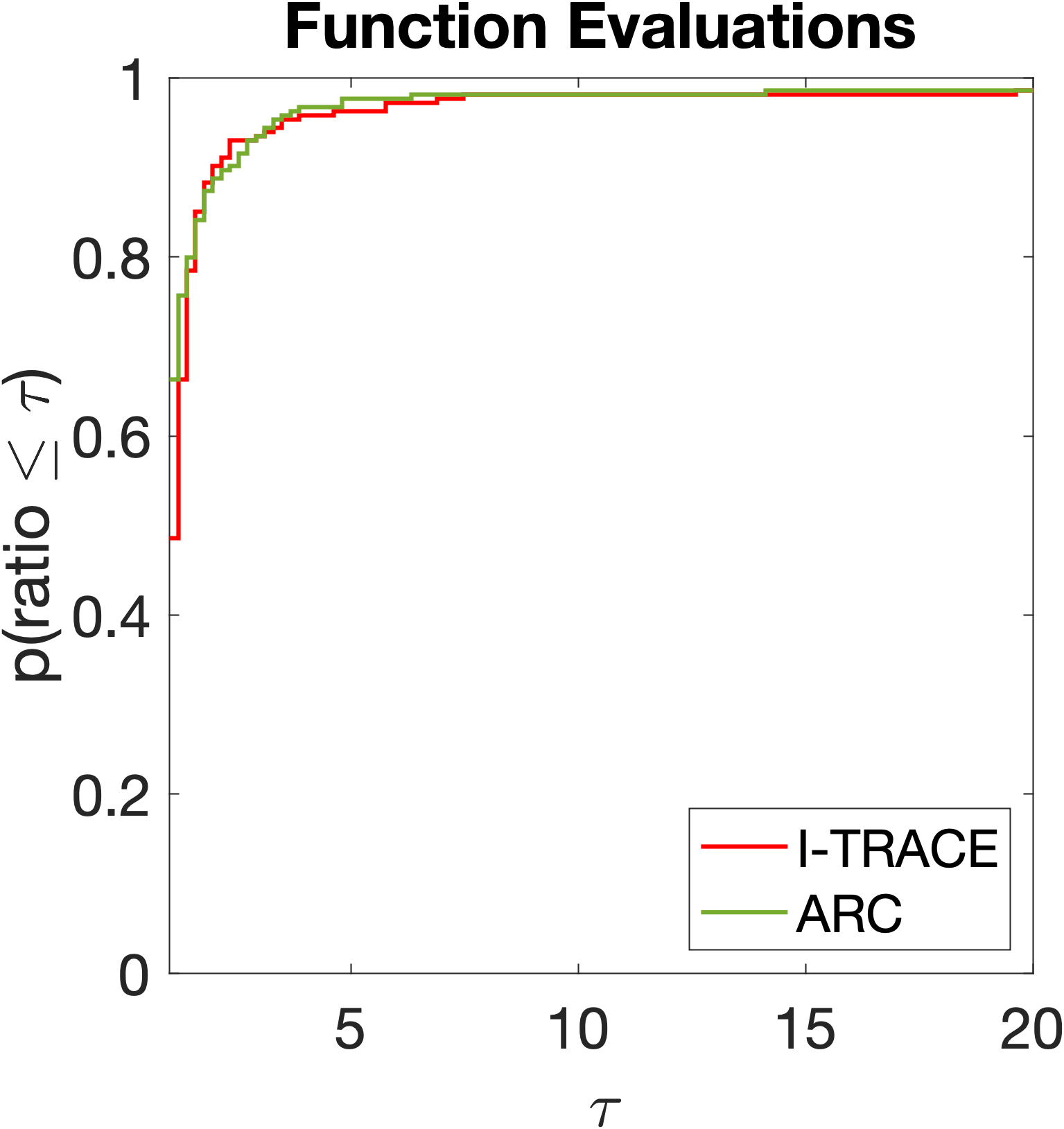}
    \includegraphics[width=5.2cm]{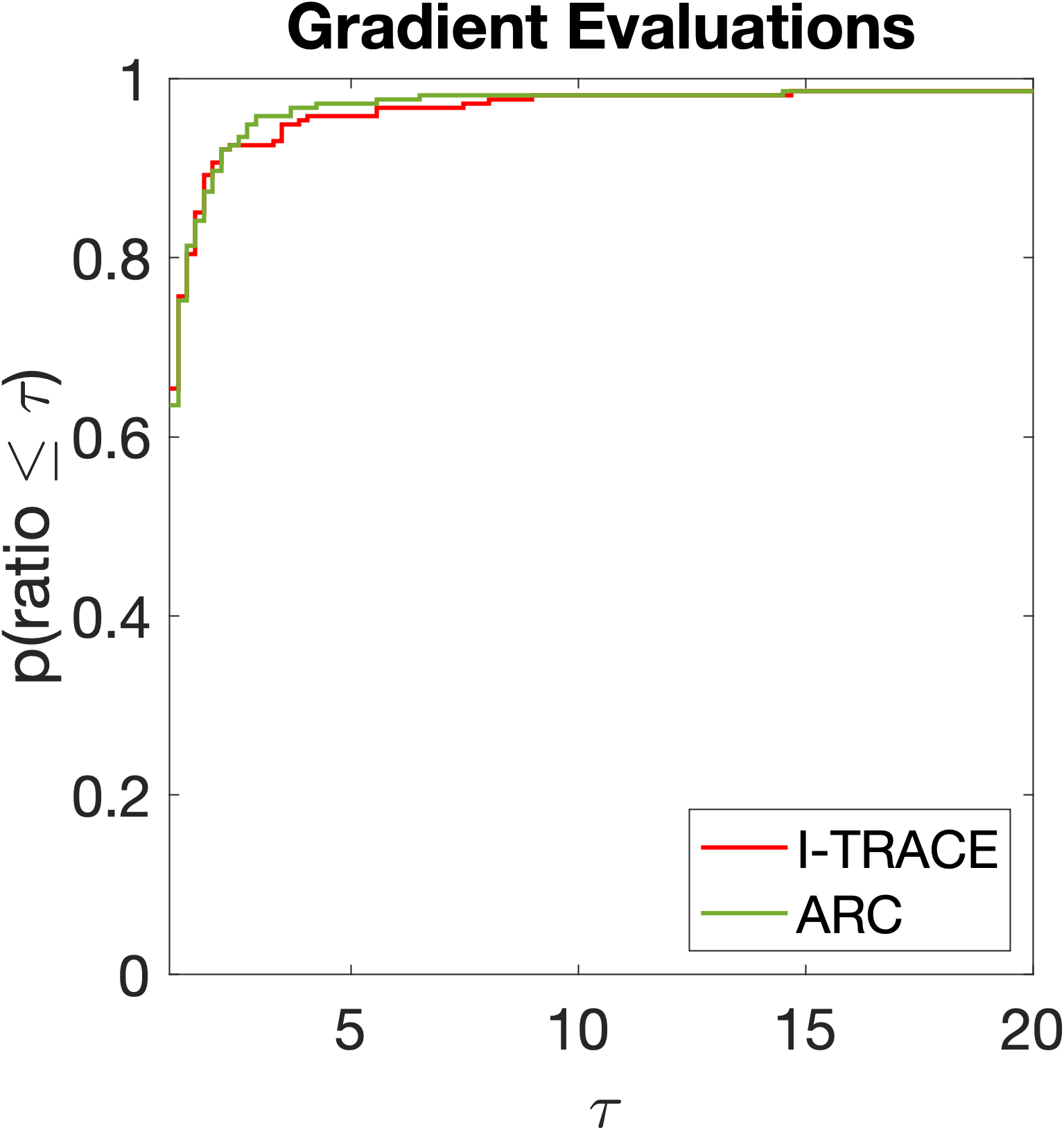}
    \includegraphics[width=5.2cm]{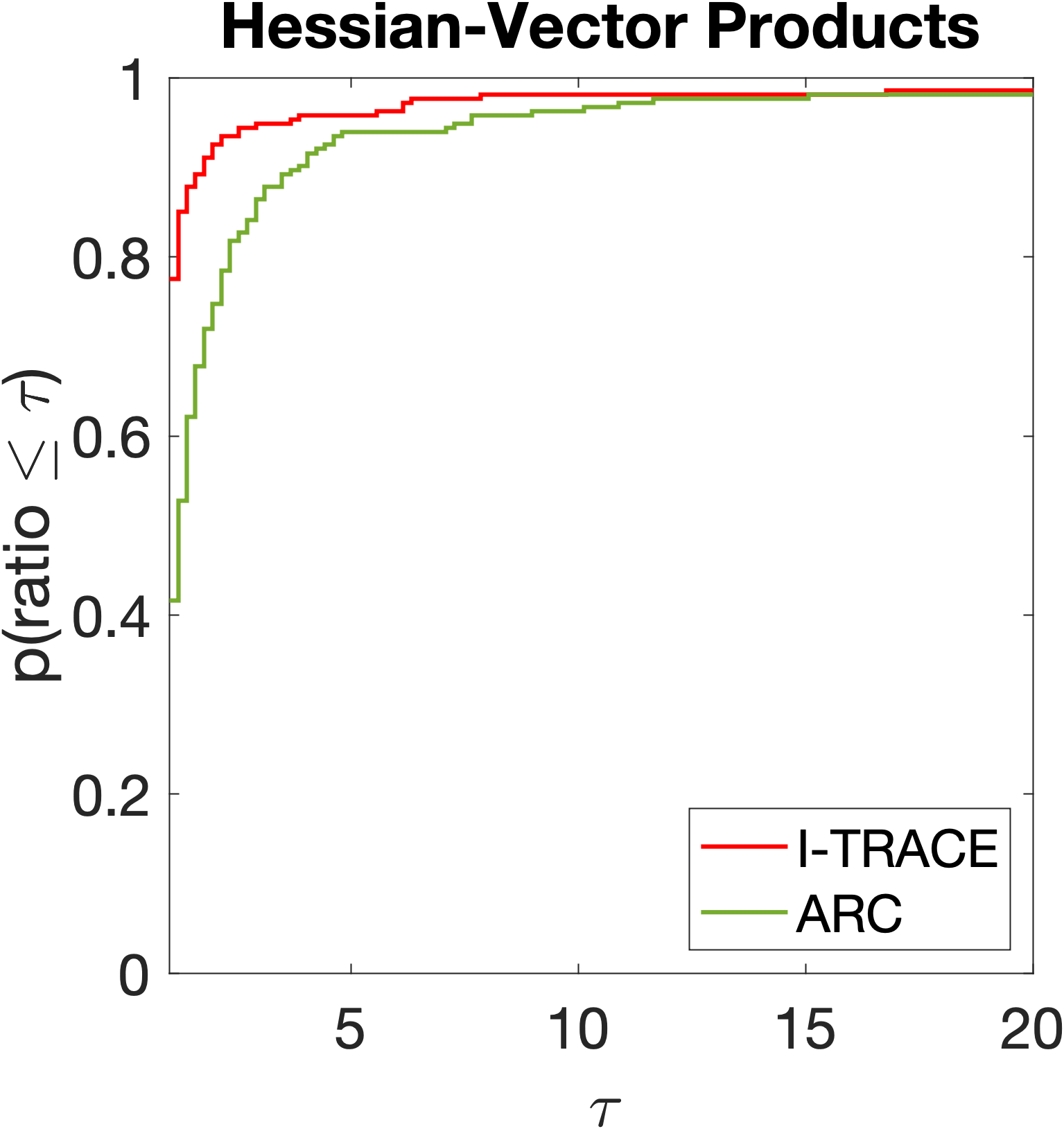}
    \caption{Performance profiles for function evaluations, gradient evaluations, and Hessian-vector products for \ITRACE{} and \ARC{} $($both with setting~3$)$ when solving 214 CUTEst problems.} \label{fig.other_pp3}
\end{figure}

We provide in Figures~\ref{fig.other_pp1}, \ref{fig.other_pp2}, and \ref{fig.other_pp3} performance profiles comparing \ITRACE{} and \ARC{} with their settings 1, 2, and 3, respectively.  Again, to focus only on relative performance for successful cases, each set of profiles only considers problems for which both algorithms were successful.  (We have already confirmed above that the reliability of the solvers were comparable for each parameter setting.)

The profiles in Figures~\ref{fig.other_pp1}, \ref{fig.other_pp2}, and \ref{fig.other_pp3} show that \ITRACE{} performs at least as well as \ARC{} across a range of parameter settings and a broad spectrum of problems.  The two algorithms perform the most alike when they both use setting~3, in which case they perform very comparably in terms of function and gradient evaluations, although \ITRACE{} performs better overall in terms of Hessian-vector products.

\section{Conclusion}\label{sec.conclusion}

We presented, analyzed, and tested a new algorithm for solving smooth unconstrained optimization problems.  The algorithm is an extension of \TRACE{}~\cite{CurtRobiSama17}, specifically one that allows the use of inexact subproblem solutions that are computed using an iterative linear algebra technique (the Lanczos algorithm, a Krylov subspace method).  The algorithm, referred to as \ITRACE{}, maintains the worst-case iteration complexity guarantees (to $\epsilon$-stationarity, as defined in \eqref{eq.eps_stationary}) and local convergence rate guarantees of \TRACE{}, but offers worst-case guarantees in terms of Hessian-vector products that can be significantly better than those offered by \TRACE{}.  Numerical experiments show that \ITRACE{} can offer better computational trade-offs than \TRACE{}, and show that \ITRACE{} is competitive with a state-of-the-art second-order method with optimal complexity guarantees to $\epsilon$-stationarity.

\bibliographystyle{plain}
\bibliography{references}

\begin{thebibliography}{10}

\bibitem{adachi2017solving}
Satoru Adachi, Satoru Iwata, Yuji Nakatsukasa, and Akiko Takeda.
\newblock {Solving the trust-region subproblem by a generalized eigenvalue
  problem}.
\newblock {\em SIAM Journal on Optimization}, 27(1):269--291, 2017.

\bibitem{BirgGardMartSantToin17}
E.~G. Birgin, J.~L. Gardenghi, J.~M. Mart{\'\i}nez, S.~A. Santos, and {\relax
  Ph}.~L. Toint.
\newblock {Worst-case evaluation complexity for unconstrained nonlinear
  optimization using high-order regularized models}.
\newblock {\em Mathematical Programming}, 163(1):359--368, 2017.

\bibitem{BirgMart17}
E.~G. Birgin and J.~M. Mart\'inez.
\newblock {The Use of Quadratic Regularization with a Cubic Descent Condition
  for Unconstrained Optimization}.
\newblock {\em SIAM Journal on Optimization}, 27(2):1049--1074, 2017.

\bibitem{carmon2018analysis}
Yair Carmon and John~C. Duchi.
\newblock {Analysis of Krylov subspace solutions of regularized nonconvex
  quadratic problems}.
\newblock In {\em Proceedings of the 32nd International Conference on Neural
  Information Processing Systems}, pages 10728--10738, 2018.

\bibitem{carmon2018accelerated}
Yair Carmon, John~C. Duchi, Oliver Hinder, and Aaron Sidford.
\newblock {Accelerated methods for nonconvex optimization}.
\newblock {\em SIAM Journal on Optimization}, 28(2):1751--1772, 2018.

\bibitem{CartGoulToin10}
C.~Cartis, N.~I.~M. Gould, and Ph.~L. Toint.
\newblock {On the Complexity of Steepest Descent, Newton's and Regularized
  Newton's Methods for Nonconvex Unconstrained Optimization Problems}.
\newblock {\em SIAM Journal on Optimization}, 20(6):2833--2852, 2010.

\bibitem{CartGoulToin11}
C.~Cartis, N.~I.~M. Gould, and {\relax Ph}.~L. Toint.
\newblock {Adaptive cubic regularisation methods for unconstrained
  optimization. Part I: Motivation, convergence and numerical results}.
\newblock {\em Mathematical Programming}, 127:245--295, 2011.

\bibitem{CartGoulToin11b}
C.~Cartis, N.~I.~M. Gould, and {\relax Ph}.~L. Toint.
\newblock {Adaptive cubic regularisation methods for unconstrained
  optimization. Part II: Worst-case function- and derivative-evaluation
  complexity}.
\newblock {\em Mathematical Programming}, 130:295--319, 2011.

\bibitem{ConGT00}
Andrew~R. Conn, Nicholas I.~M. Gould, and {\relax Ph}ilippe~L. Toint.
\newblock {A primal-dual algorithm for minimizing a non-convex function subject
  to bound and linear equality constraints}.
\newblock In {\em Nonlinear optimization and related topics (Erice, 1998)},
  volume~36 of {\em Appl. Optim.}, pages 15--49. Kluwer Acad. Publ., Dordrecht,
  2000.

\bibitem{ConGT00a}
Andrew~R. Conn, Nicholas I.~M. Gould, and {\relax Ph}ilippe~L. Toint.
\newblock {\em {Trust-Region Methods}}.
\newblock Society for Industrial and Applied Mathematics (SIAM), Philadelphia,
  PA, 2000.

\bibitem{CurtLubbRobi18}
F.~E. Curtis, Z.~Lubberts, and D.~P. Robinson.
\newblock {Concise complexity analyses for trust region methods}.
\newblock {\em Optimization Letters}, 2018.

\bibitem{CurtRobiSama17}
F.~E. Curtis, D.~P. Robinson, and M.~Samadi.
\newblock {A trust region algorithm with a worst-case iteration complexity of
  $\mathcal{O}(\epsilon^{-3/2})$ for nonconvex optimization}.
\newblock {\em Mathematical Programming}, 162(1):1--32, 2017.

\bibitem{CurtRobiSama18}
F.~E. Curtis, D.~P. Robinson, and M.~Samadi.
\newblock {An inexact regularized Newton framework with a worst-case iteration
  complexity of $\mathcal{O}(\epsilon^{-3/2})$ for nonconvex optimization}.
\newblock {\em IMA Journal on Numerical Analysis}, 2018.

\bibitem{curtis2021trust}
Frank~E. Curtis, Daniel~P. Robinson, Cl{\'e}ment~W. Royer, and Stephen~J.
  Wright.
\newblock {Trust-region Newton-CG with strong second-order complexity
  guarantees for nonconvex optimization}.
\newblock {\em SIAM Journal on Optimization}, 31(1):518--544, 2021.

\bibitem{DemES82}
Ron~S. Dembo, Stanley~C. Eisenstat, and Trond Steihaug.
\newblock {Inexact Newton methods}.
\newblock {\em SIAM J. Numer. Anal.}, 19(2):400--408, 1982.

\bibitem{DolM02}
Elizabeth~D. Dolan and Jorge~J. Mor{\'e}.
\newblock {Benchmarking optimization software with performance profiles}.
\newblock {\em Math. Program.}, 91(2, Ser. A):201--213, 2002.

\bibitem{Duss17}
J.-P. Dussault.
\newblock {ARCq: a new adaptive regularization by cubics}.
\newblock {\em Optimization Methods and Software},
  doi:10.1080/10556788.2017.1322080, 2017.

\bibitem{DussOrba15}
J.-P. Dussault and D.~Orban.
\newblock {Scalable adaptive cubic regularization methods}.
\newblock Technical Report G-2015-109, GERAD, 2015.

\bibitem{GoulPorcToin12}
N.~I.~M. Gould, M.~Porcelli, and P.~L. Toint.
\newblock {Updating the Regularization Parameter in the Adaptive Cubic
  Regularization Algorithm}.
\newblock {\em Computational Optimization and Applications}, 53(1):1--22, Sep
  2012.

\bibitem{GouLRT99}
Nicholas I.~M. Gould, Stefano Lucidi, Massimo Roma, and {\relax Ph}ilippe~L.
  Toint.
\newblock {Solving the trust-region subproblem using the Lanczos method}.
\newblock {\em SIAM J. Optim.}, 9(2):504--525, 1999.

\bibitem{GouOT13}
Nicholas I.~M. Gould, D.~Orban, and {\relax Ph}ilippe~L. Toint.
\newblock {CUTEst: a Constrained and Unconstrained Testing Environment with
  safe threads}.
\newblock Technical report, Rutherford Appleton Laboratory, Chilton, England,
  2013.

\bibitem{gould2020error}
Nicholas I.~M. Gould and Valeria Simoncini.
\newblock {Error estimates for iterative algorithms for minimizing regularized
  quadratic subproblems}.
\newblock {\em Optimization Methods and Software}, 35(2):304--328, 2020.

\bibitem{GrapYuanYuan15}
G.~N. Grapiglia, J.~Yuan, and Y.~Yuan.
\newblock {On the convergence and worst-case complexity of trust-region and
  regularization methods for unconstrained optimization}.
\newblock {\em Mathematical Programming}, 152(1-2):491--520, 2015.

\bibitem{Grie81}
A.~Griewank.
\newblock {The Modification of Newton's Method for Unconstrained Optimization
  by Bounding Cubic Terms}.
\newblock Technical Report NA/12, Department of Applied Mathematics and
  Theoretical Physics, University of Cambridge, 1981.

\bibitem{JiaW21}
Zhongxiao Jia and Fa~Wang.
\newblock {The convergence of the generalized Lanczos trust-region method for
  the trust-region subproblem}.
\newblock {\em SIAM Journal on Optimization}, 31(1):887--914, 2021.

\bibitem{kuczynski1994probabilistic}
J~Kuczy{\'n}ski and H~Wo{\'z}niakowski.
\newblock {Probabilistic bounds on the extremal eigenvalues and condition
  number by the Lanczos algorithm}.
\newblock {\em SIAM Journal on Matrix Analysis and Applications},
  15(2):672--691, 1994.

\bibitem{kuczynski1992estimating}
Jacek Kuczy{\'n}ski and Henryk Wo{\'z}niakowski.
\newblock {Estimating the largest eigenvalue by the power and Lanczos
  algorithms with a random start}.
\newblock {\em SIAM journal on matrix analysis and applications},
  13(4):1094--1122, 1992.

\bibitem{MorS83}
Jorge~J. Mor\'{e} and Danny~C. Sorensen.
\newblock {Computing a trust region step}.
\newblock {\em SIAM J. Sci. and Statist. Comput.}, 4:553--572, 1983.

\bibitem{NestPoly06}
{Yu.} Nesterov and B.~T. Polyak.
\newblock {Cubic regularization of Newton's method and its global performance}.
\newblock {\em Mathematical Programming}, 108(1):117--205, 2006.

\bibitem{NoceWrig06}
J.~Nocedal and S.~J. Wright.
\newblock {\em {Numerical Optimization}}.
\newblock Springer Series in Operations Research. Springer, New York, NY, USA,
  second edition, 2006.

\bibitem{royer2020newton}
Cl{\'e}ment~W. Royer, Michael O’Neill, and Stephen~J. Wright.
\newblock {A Newton-CG algorithm with complexity guarantees for smooth
  unconstrained optimization}.
\newblock {\em Mathematical Programming}, 180(1):451--488, 2020.

\bibitem{royer2018complexity}
Cl{\'e}ment~W. Royer and Stephen~J. Wright.
\newblock {Complexity analysis of second-order line-search algorithms for
  smooth nonconvex optimization}.
\newblock {\em SIAM Journal on Optimization}, 28(2):1448--1477, 2018.

\bibitem{Ste83}
Trond Steihaug.
\newblock {The conjugate gradient method and trust regions in large scale
  optimization}.
\newblock {\em SIAM J. Numer. Anal.}, 20:626--637, 1983.

\bibitem{Toi81}
{\relax Ph}ilippe~L. Toint.
\newblock {Towards an efficient sparsity exploiting Newton method for
  minimization}.
\newblock In I.~S. Duff, editor, {\em Sparse Matrices and Their Uses}, pages
  57--88, London and New York, 1981. Academic Press.

\bibitem{ZhaSL17}
Lei-Hong Zhang, Chungen Shen, and Ren-Cang Li.
\newblock {On the generalized Lanczos trust-region method}.
\newblock {\em SIAM Journal on Optimization}, 27(3):2110--2142, 2017.

\end{thebibliography}

\end{document}